\documentclass[11pt]{article}

\usepackage[
backend=biber,
style=numeric-comp,
sorting=nty
]{biblatex}
\usepackage[margin=1in]{geometry}
\usepackage[utf8]{inputenc}
\usepackage{csquotes}
\usepackage{amsfonts,amsthm,amsmath,amssymb,accents,mathtools,physics} 
\usepackage{mathrsfs,amsbsy} 
\usepackage{enumitem}
\usepackage{subfig}
\usepackage{graphicx}
\usepackage{float}
\usepackage{xcolor}
\usepackage{algorithm}
\usepackage[noend]{algpseudocode}
\usepackage{hyperref}
\usepackage{authblk} 

\usepackage{textcomp} 
\usepackage{algorithm}
\usepackage{algorithmicx}
\usepackage{relsize}

\usepackage{bbm}
\usepackage{dsfont}

\numberwithin{equation}{section}

\usepackage{yhmath}
 \usepackage{booktabs} 

\newtheorem{rem}{Remark}[section]

\theoremstyle{definition}

\newcommand{\vperp}{v_{\perp}}
\newcommand{\vpar}{v_{||}}
\newcommand{\upar}{u_{||}}
\newcommand{\Ta}{T_{\alpha}}
\newcommand{\Te}{T_e}

\newcommand{\bfV}{\mathbf{V}}

\newcommand{\bfK}{\mathbf{K}}

\usepackage{color}

\hypersetup{
    colorlinks=true,
    linkcolor=blue,
    citecolor=blue,
    filecolor=blue,      
    urlcolor=blue,
}

\title{An adaptive and conservative low-rank IMEX solver for the hybrid ion Vlasov-Fokker-Planck and fluid electron system}

\author[1]{Joseph Nakao}
\author[2]{Dylan T. Jacobs}
\author[3]{William Taitano}

\affil[1]{Department of Mathematics and Statistics, Swarthmore College, Swarthmore, PA, USA\newline
Email: jnakao1@swarthmore.edu, ORCID iD: 0009-0008-6589-4013}
\affil[2]{Department of Mathematics and Statistics, Swarthmore College, Swarthmore, PA, USA\newline
Email: djacobs2@swarthmore.edu, ORCID iD: 0009-0009-8887-2703}
\affil[3]{Theoretical Division, Los Alamos National Laboratory, Los Alamos, NM, USA\newline
Email: taitano@lanl.gov, ORCID iD: 0000-0002-2369-0935}

\date{}

\begin{document}

\maketitle
\let\thefootnote\relax\footnotetext{\hspace*{-1.8em}\textbf{Corresponding author:} Joseph Nakao\\
\textbf{MSC2020 Numbers:} 65M06\\
\textbf{Keywords:} Vlasov-Fokker-Planck, Chang-Cooper, implicit-explicit, rank-adaptive, conservative SVD
}


\begin{abstract}
\noindent We present a macroscopically conservative, rank-adaptive method for solving a hybrid Vlasov-Fokker-Planck (VFP) equation in which the ions are treated kinetically, and the electrons are treated as a fluid. Solving this system poses several coupled computational challenges: the curse of dimensionality, conservation of macroscopic quantities, stiffness arising from the multiscale nature of the model, and structure preservation. To address these challenges, we combine several established methods, each targeting a specific difficulty, into a unified framework for solving this nonlinear system. Specifically, we employ the recent Reduced Augmentation Implicit Low-rank (RAIL) method for rank adaptivity in velocity space to combat the curse of dimensionality, high-order implicit-explicit time stepping to handle the multiscale stiffness, and the Local Macroscopic Conservative (LoMaC) procedure for conservative truncation of the solution. The RAIL and LoMaC methods are extended from Cartesian to cylindrical coordinates. The resulting framework accommodates implicit rank-adaptive time integration while remaining conservative and leveraging structure-preserving discretizations. We verify the scheme on a suite of test problems and simulate a standing shock to demonstrate the importance of conserving the macroscopic quantities.
\end{abstract}


\section{Introduction}
The Vlasov-Landau-Rosenbluth-Fokker-Planck (VLRFP) equation describes the evolution of the plasma particle distribution function in six-dimensional position-velocity phase space subject to Coulomb collisions and electromagnetic fields. Coupled to Maxwell's equations, it is regarded as the first-principles model for weakly coupled plasmas across all collisionality regimes, with applications ranging from fusion energy and high energy density physics to space weather and semiconductor manufacturing. The equation is integro-differential, multiscale in both time and phase space, and preserves important symmetries, including conservation of mass, momentum, and energy, as well as the Boltzmann $\mathcal{H}$-theorem. Its high-dimensional nature further incurs the so-called \emph{curse of dimensionality}, whereby computational and storage costs grow exponentially with dimension. These properties have made the equation a focus of intense study in both the plasma physics and applied mathematics communities. The resulting literature is far too extensive to review comprehensively here; instead, we survey a representative subset to provide sufficient context for our proposed approach.

Techniques for preserving discrete conservation include approaches based on projection and on nonlinear constraints. In both, the continuum symmetries necessary for conservation are first identified, and the discretization truncation errors that violate these symmetries are then either projected out directly at the solution level through post-processing (projection methods) or eliminated by modifying the discrete operators (nonlinear constraint methods). Projection-based methods include the Local Macroscopic Conservative (LoMaC) method, which relies on a dual representation of the solution wherein the underlying kinetic equation is co-evolved with self-consistent moment equations. Although the kinetic equation may not directly preserve discrete conservation symmetries, the moment equations, being written in divergence form, do. The closures for viscous stress and heat flux follow from taking moments of the kinetic advection operator, which likewise assumes a divergence form in space. The kinetic solution is then coupled to the moment equations through a post-processing operation per time-step. The nonlinear constraint approach, by contrast, does not require an accompanying moment equation. Instead, this approach identifies the underlying discrete symmetries of the operator and uses backward error analysis to quantify the truncation error introduced by the numerical method, then modifies the operator accordingly. Alternative structure-preserving techniques rely on structural transformation of the underlying kinetic equation, which separates the conserved moments from the evolution of the distribution function and isolates conservation into the associated moment equation. The micro-macro decomposition (MMD) is the most well-known of these methods \cite{liu_cmp_2004_mmd, lemou_jsc_2008_mmd, bennoune_jcp_2008_mmd, gamba_jcp_2019_mmd}, wherein the perturbation to the distribution function is evolved without carrying information on mass, momentum, or energy. More recently, the conditional kinetics (CK) approach has been introduced, which relies on a coordinate transformation in the velocity space to separate out the conserved moments from the evolution of the conditional distribution function \cite{taitano_jpp_2026_ck}. In both approaches, associated moment equations are evolved. While both approaches share similarities with LoMaC, they differ in a key aspect: whereas LoMaC uses moment equations as an auxiliary system to inform a downstream projection operator for conservation, MMD and CK self-consistently evolve the moment equations as part of the dynamical system. 

To address the multiscale challenges of the Vlasov-Fokker-Planck equations in both space and time, various frameworks have been devised, including MMD and CK. In both approaches, stiff time scales associated with collisions and collective phenomena such as plasma waves are segregated from the kinetic equations and isolated into lower-dimensional, more tractable moment equations. Alternatively, the unified gas kinetic scheme (UGKS) \cite{xu2010unified,Chen2015} has been introduced, which uses the exact characteristics-based solution of the Boltzmann transport equation to inform the numerical flux, ensuring a smooth transition between the free-streaming and Navier-Stokes limits. On the other hand, higher-order low-order (HOLO) methods \cite{willert_neutraon_holo_jcp_2014, taitano_bgk_holo_jctt_2014, taitano_vfp_holo_jcp_2015, holo_review_2017} introduce moment equations as auxiliary systems to capture fast collective dynamics, thereby relaxing the need to invert otherwise dense operators (e.g., collision integrals and moment coupling with Maxwell's equations). Once the auxiliary moment equations are solved, the stiff operators in the underlying kinetic equation are parameterized using the auxiliary moment solution, sparing the kinetic solver from inverting these dense operators directly. Both the HOLO and UGKS methods have been successfully used in integrated simulations of complex engineering systems including modeling of inertial confinement fusion (ICF) systems \cite{mannion_pre_2023_omega, johnson_pre_2024_omega} and high-altitude reentry vehicles \cite{Wei2026}. To address multiscale challenges in space, particularly recovering asymptotic reduced-order limiting forms such as the Euler/Navier-Stokes equations and ambipolarity conditions, MMD and CK methods analytically decouple the conserved moments from the evolution of the distribution function, thereby isolating critical symmetries into the moment-field subsystem and simplifying the work required of the discrete scheme in the high-dimensional kinetic equation. Beyond the MMD and CK approaches, significant work has been performed on directly recovering appropriate asymptotic limiting solutions by discretizing the original Vlasov-Fokker-Planck equation. These include the Chang-Cooper scheme \cite{ChangC1970practical, PareschiZ2018structure} and the quadrature-corrected moment (QCM) scheme \cite{TaitanoCS2017equilibrium}, which ensure that the exact equilibrium is recovered in the discrete formulation by appropriately identifying the correct analytical kernel of the collision operator.

Furthermore, considerable work has been devoted to reducing the computational complexity of kinetic equations. Ref. \cite{kolobov_aip_cp_2019_amr_vlasov} has explored classic adaptive-mesh-refinement (AMR) techniques for this purpose, while Refs. \cite{taitano_2016_jcp_r_adaptivity,taitano_2018_jcp_r_adaptivity,taitano_2021_jcp_r_adaptivity} have developed a physics-based r-adaptive phase-space grid, wherein the velocity-space grid is dynamically adapted based on the local, instantaneous temperature and drift velocity to track important bulk features of the distribution. Alternatively, complexity-reduction techniques based on low-rank factorization of the kinetic solution have gained considerable attention over the past decade. Recent advances in matrix factorization and, more generally, low-rank tensor decompositions \cite{hackbusch2009new,grasedyck2013literature,kormann2017low,kolda2009tensor,bachmayr2023low} have been increasingly applied to the curse of dimensionality inherent in kinetic equations. By factorizing the underlying high-dimensional object into factor matrices and compressed core/coefficient tensors, the storage complexity is sharply reduced. Within the context of \textit{purely kinetic} Fokker-Planck and Vlasov-Fokker-Planck equations, many low-rank methods have recently been proposed \cite{wang2026implicit,wang2026dynamical,zhang2026separable,coughlin2024robust,coughlin2022efficient,KahzaTQC2024krylov,NakaoQE2025reduced,chertkov2021solution,dolgov2012fast,dektor2021rank,carrel2026sketch,nelson2026local}. Ref. \cite{einkemmer2025review} provides a comprehensive review of recent low-rank methods used for kinetic simulations.

In this manuscript, we present a multiscale, conservative, and discrete equilibrium-preserving low-rank algorithm for solving the one-dimension-in-space, two-dimensions-in-cylindrical-velocity (1D2V) \emph{hybrid} system of kinetic ions and fluid electrons. Under asymptotically limiting conditions where the plasma Debye length vanishes relative to local gradient length scales, the Vlasov-Amp{\`e}re equations for electrostatic plasmas with irrotational electric fields can be simplified by expressing the field algebraically in terms of moments of the distribution function via Ohm's law. In this limit, electron densities and drift velocities can also be represented algebraically in terms of ion moments through the quasi-neutrality and ambipolarity relations, allowing one to analytically remove ultrafast electron advection and plasma wave time scales from the system, only requiring the energy equation to be evolved for the electrons (instead of the full Vlasov-Fokker-Planck). This so-called \emph{hybrid} formulation is more tractable for large-scale problems and is a popular choice among practitioners. For the ion-ion and ion-electron Fokker-Planck collision operator, we adopt a particular reduced form often referred to in the literature as the Lenard-Bernstein form, which assumes a simpler scalar isotropic diffusion tensor and a friction coefficient that depends linearly on the velocity coordinate. This simplified form is often discussed alongside the general Landau-Rosenbluth form, which involves solving coupled Poisson equations for the Rosenbluth potentials \cite{rosenbluth_1957_prl_fp, taitano2015mass} and computing their Hessian and gradient to evaluate the collisional coefficients. Despite its simplification, the operator retains the essential physical structure: for like-species collisions, with the drift and temperature set to the self-consistent moments of the distribution, it conserves mass, momentum, and energy and satisfies Boltzmann's $\mathcal{H}$-theorem, while the ion-electron collisions conserve mass and exchange energy with the electrons so that these quantities are conserved for the coupled ion-electron system as a whole. We exploit the time-scale separation between the fast collision operator, which acts in the velocity subspace, and the slower ion advection scale in physical space to devise a high-order implicit-explicit (IMEX) Runge-Kutta (RK) integrator for the coupled system. We invoke a low-rank ansatz of the distribution function in velocity space, employing the singular value decomposition for compression. In this representation, each physical-space grid point stores the factorized matrices that represent the velocity distribution. This approach mitigates a common pitfall of global low-rank methods, in which a distribution that is not separable on the global velocity grid---for instance a Maxwellian whose drift and temperature vary in space, or the cylindrical $(\vperp,\vpar)$ representation of a drifting Maxwellian---can require a high rank to represent. Storing the factorization per spatial node localizes the drift so that each node sees a nearly centered, low-rank distribution. Furthermore, to treat the stiff velocity-space operators implicitly with high-order accuracy in time, we extend the recent Reduced Augmentation Implicit Low-rank (RAIL) solver to the coupled system in cylindrical velocity coordinates. Separately, the fluid electron equation supports a stiff thermal conduction term, which is nonlinear, and is therefore treated implicitly using a Newton solver. Finally, we ensure conservation by employing a LoMaC scheme generalized to cylindrical coordinates, wherein the ion moment solution is solved in tandem with the electron energy equation. While the RAIL method is used as the implicit integrator, our proposed method can accommodate other similar low-rank implicit integrators.

The remainder of the manuscript is organized as follows. Section 2 presents the coupled ion Vlasov-Fokker-Planck and fluid electron equations under consideration. Section 3 presents our numerical scheme, covering the individual low-rank kinetic ion solver, conservative projectors, time integrators, and the coupled fluid solver for the ion moment and electron energy equations. Section 4 presents the numerical results, and we conclude in Section 5.


\section{The model of interest}
We consider the Vlasov-Fokker-Planck (VFP) equation describing the dynamics of weakly coupled collisional plasmas. The VFP equation is given by
\begin{equation}
    \frac{\partial f_{\alpha}}{\partial t} + \mathbf{v} \cdot \nabla_x f_{\alpha} + \frac{q_{\alpha}}{m_{\alpha}}(\mathbf{E}+\mathbf{v}\times\mathbf{B}) \cdot \nabla_v f_{\alpha} = \sum_{\beta=1}^{N_s}C_{\alpha \beta}
    \label{eq:vfp-full-model}
\end{equation}
where $f_{\alpha}(\mathbf{x},\mathbf{v},t)$ is the particle distribution function for particle species $\alpha$ in physical space $\mathbf{x}\in\Omega_x\subset\mathbb{R}^3$, velocity space $\mathbf{v}\in\mathbb{R}^3$, and time $t\in\mathbb{R}_+$. $\mathbf{E}(\mathbf{x},t)$ is the electric field, $\mathbf{B}(\mathbf{x},t)$ is the magnetic field, $q_{\alpha}$ is the particle charge, $m_{\alpha}$ is the particle mass, and $C_{\alpha\beta}$ is the Lenard-Bernstein-Fokker-Planck (LBFP) collision operator describing the small-angle binary collisions between species $\alpha$ and species $\beta$ \cite{LenardB1958},
\begin{equation}\label{eq:LBFP}
    C_{\alpha\beta} = \nu_{\alpha\beta}\nabla_v \cdot \left(\frac{T_{\beta}}{m_{\alpha}}\nabla_v f_{\alpha} + (\mathbf{v} - \mathbf{u}_{\beta})f_{\alpha}\right),
\end{equation}
where constant $\nu$ is the collision frequency, $\mathbf{u}$ is the drift/bulk velocity, and $T$ is the temperature. The LBFP operator \eqref{eq:LBFP} is also known as the Dougherty-Fokker-Planck operator \cite{Dougherty1964}. Although the LBFP operator is linear with respect to the independent variables, the macroscopic quantities are nonlinear functionals of the distribution function. When the drift $\mathbf{u}_{\beta}$ and temperature $T_{\beta}$ in the like-species operator $C_{\alpha\alpha}$ are set to the self-consistent moments of $f_{\alpha}$, the operator conserves mass, momentum, and energy and satisfies the $\mathcal{H}$-theorem. The inter-species operator $C_{\alpha e}$ conserves mass but exchanges energy with the electron fluid through the term $W_{e\alpha}$ defined in equation \eqref{eq:fluidelectron}; the total mass, momentum, and energy of the coupled ion-electron system are conserved. The number density, drift velocity, and temperature are defined by the first few moments of the distribution function,
\begin{subequations}
\begin{equation}
    n(\mathbf{x},t) = \big\langle f(\mathbf{x},\mathbf{v},t),1\big\rangle_v,
\end{equation}
\begin{equation}
    \mathbf{u}(\mathbf{x},t) = \frac{1}{n}\big\langle f(\mathbf{x},\mathbf{v},t),\mathbf{v}\big\rangle_v,
\end{equation}
\begin{equation}
    T(\mathbf{x},t) = \frac{m}{3n}\big\langle f(\mathbf{x},\mathbf{v},t),|\mathbf{v}-\mathbf{u}|^2\big\rangle_v,
\end{equation}
\end{subequations}
where $\langle\cdot,\cdot\rangle_v$ is the unweighted $L^2$ inner product over $\mathbb{R}^3$ velocity space. The factor of $3$ in the temperature reflects the three translational degrees of freedom; in the azimuthally symmetric $(\vperp,\vpar)$ system the perpendicular direction carries two and the parallel direction one.

In this study, we consider 1D-2V phase space $(x,\vperp,\vpar)$ under a 1D planar geometry in physical space, and a 2V cylindrically symmetric coordinate system in velocity space, where $\mathbf{E}=E_{||}\hat{x}$ and $\hat{\vpar}$ point in the same direction. The inner product under this 2V coordinate system is now defined as $\langle A,B\rangle_v\coloneqq 2\pi\int_{-\infty}^{\infty}\int_0^{\infty}{AB\vperp d\vperp d\vpar}$. As an initial step towards solving the general VFP equation in our proposed framework, we assume zero magnetic field and a single ion species $\alpha$.

Following the work of \cite{taitano_2018_jcp_r_adaptivity,TaitanoCS2021conservative}, we keep the kinetic model for the ions, but use a fluid model for the electron pressure $p_e=n_eT_e$ and assume parallel drift $u_{\perp}=0$ so that $\mathbf{u}=u_{||}\hat{x}$. Furthermore, we assume quasi-neutrality $n=n_{\alpha}=n_e$, ambipolarity $\mathbf{u}=\mathbf{u}_{\alpha}=\mathbf{u}_e$, and use Ohm's law to express the electric field $E_{||}=\frac{1}{q_e n_e}\frac{\partial p_e}{\partial x}$. The resulting \textit{hybrid kinetic-ion, fluid-electron model} describing the distribution of ion species $\alpha$, denoted by $f\equiv f_{\alpha}$, is
\begin{equation}\label{eq:VFP1d2v}
    \frac{\partial f}{\partial t} + \vpar\frac{\partial f}{\partial x} + \frac{q_{\alpha}}{m_{\alpha}}E_{||}\frac{\partial f}{\partial\vpar} = C_{\alpha\alpha} + C_{\alpha e},
\end{equation}
\begin{equation}\label{eq:fluidelectron}
    \frac{3}{2}\frac{\partial p_e}{\partial t} + \frac{5}{2}\frac{\partial}{\partial x}(u_{||}p_e) - u_{||}\frac{\partial p_e}{\partial x} - \frac{\partial}{\partial x}\Big(\kappa_{||}\frac{\partial T_e}{\partial x}\Big) = W_{e\alpha},
\end{equation}
where $C_{\alpha\alpha}$ and $C_{\alpha e}$ are defined by equation \eqref{eq:LBFP}, $W_{e\alpha} = -\left\langle\frac{m_{\alpha}|\mathbf{v}|^2}{2},C_{\alpha e}\right\rangle_v = 3\nu_{\alpha e}n(T_{\alpha}-T_e)$ describes the electron-ion energy exchange, and $\kappa_{||}$ is the electron thermal conductivity. The formulas for the collision frequencies and thermal conductivity can be found in \cite{HazeltineM2003plasma,taitano_2018_jcp_r_adaptivity,TaitanoCS2021conservative,SimakovM2014electron}, as
\begin{equation}\label{eq:collisionfreq}
    \nu_{\alpha\alpha} = \frac{\zeta n}{\sqrt{m_{\alpha}}T_{\alpha}^{3/2}},\qquad\qquad \nu_{\alpha e} = \frac{\zeta\sqrt{2m_e}n}{m_{\alpha}T_e^{3/2}},
\end{equation}
\begin{equation}\label{eq:thermalcond}
    \kappa_{||} = \frac{3.2T_e^{5/2}}{\zeta\sqrt{2m_e}},
\end{equation}
for some fixed constant $\zeta$ that absorbs the Coulomb logarithm and the remaining collisional constants, and where the numerical coefficient $3.2$ in \eqref{eq:thermalcond} is the classical Braginskii/Spitzer-H\"arm parallel electron thermal-conduction coefficient \cite{HazeltineM2003plasma,SimakovM2014electron}. We introduce the total specific energy $U$ (kinetic energy per unit mass), so that $nU\coloneqq\left\langle f,\frac{\norm{\mathbf{v}}^2}{2}\right\rangle_v$ is the kinetic energy density. Using the expansion $\norm{\mathbf{v}}^2=(v_{||}-u_{||})^2+u_{||}^2+v_{\perp}^2+2u_{||}(v_{||}-u_{||})$, the cross term vanishes because $\langle f,v_{||}-u_{||}\rangle_v=0$ by the definition of $u_{||}$, and $\langle f,(v_{||}-u_{||})^2+v_{\perp}^2\rangle_v=3nT_{\alpha}/m_{\alpha}$ by the temperature definition; this yields the relationship
\begin{equation}\label{eq:UTrelationship}
    2nU = \frac{3nT_{\alpha}}{m_{\alpha}} + nu_{||}^2.
\end{equation}

We nondimensionalize equations \eqref{eq:VFP1d2v}-\eqref{eq:fluidelectron} using the following arbitrary reference quantities:
\begin{align}
\begin{split}
    n^*,\qquad T^*,\qquad m^*,\qquad &q^*,\qquad u^*=\sqrt{\frac{T^*}{m^*}},\\
    t^*=(\nu^*)^{-1}=\frac{\sqrt{m^*}(T^*)^{3/2}}{\zeta n^*},\qquad &L^*=u^*t^*,\qquad f^*=\frac{n^*}{(u^*)^3},\\
    E^*=\frac{T^*}{q^*L^*},\qquad& \kappa^*=\frac{(T^*)^{5/2}}{\zeta\sqrt{m^*}}.
\end{split}
\end{align}

The resulting nondimensionalized equations are identical to the equations already presented, with the exception of equations \eqref{eq:collisionfreq}-\eqref{eq:thermalcond} for which the constant $\zeta$ is scaled out. We scale the ion mass and ion charge to unity, setting $m^*=m_{\alpha}$ and $q^*=q_{\alpha}$. Note that $m_{e}/m_{\alpha}=1/1836$, and the nondimensionalized electron charge is $q_e=-1$. The remainder of this paper assumes the nondimensionalized equations and quantities.

Taking the first few moments of equation \eqref{eq:VFP1d2v}, adding the electron energy equation \eqref{eq:fluidelectron} to the ion energy equation, and integrating over $\Omega_x$, the balance equations for the total mass, momentum, and energy are respectively
\begin{subequations}\label{eq:ionfluidsystem}
\begin{align}
    0 &= \frac{d}{dt}\int_{\Omega_x}{n(x,t)dx} + \Big[nu_{||}\Big]_{\partial\Omega_x},\\
    0 &= \frac{d}{dt}\int_{\Omega_x}{(nu_{||})(x,t)dx} + \Bigg[S_{\alpha}+\frac{nT_e}{m_{\alpha}}\Bigg]_{\partial\Omega_x}\\
    0 &= \frac{d}{dt}\int_{\Omega_x}{\Big(nU+\frac{3}{2}nT_e\Big)(x,t)dx} + \Bigg[Q_{\alpha}+\frac{5}{2}u_{||}nT_e - \kappa_{||}\frac{\partial T_e}{\partial x}\Bigg]_{\partial\Omega_x},
\end{align}
\end{subequations}
where the kinetic momentum flux $S_{\alpha}$ and energy flux $Q_{\alpha}$ from the Vlasov term are
\begin{subequations}\label{eq:kineticfluxes}
\begin{equation}
    S_{\alpha} = \Big\langle v_{||}f(x,v_{\perp},v_{||},t),v_{||}\Big\rangle_v,
\end{equation}
\begin{equation}
    Q_{\alpha} = \Bigg\langle v_{||}f(x,v_{\perp},v_{||},t),\frac{v_{\perp}^2 + v_{||}^2}{2}\Bigg\rangle_v.
\end{equation}
\end{subequations}

\begin{rem}
We verify the sign convention of the electric field force. With the nondimensional charges $q_{\alpha}=1$ and $q_e=-1$, Ohm's law gives $E_{||}=\frac{1}{q_en_e}\frac{\partial p_e}{\partial x}=-\frac{1}{n}\frac{\partial(nT_e)}{\partial x}$. The $v_{||}$-moment of the acceleration term $\frac{q_{\alpha}}{m_{\alpha}}E_{||}\frac{\partial f}{\partial\vpar}$ contributes $\frac{q_{\alpha}}{m_{\alpha}}nE_{||}=-\frac{\partial}{\partial x}\!\left(\frac{nT_e}{m_{\alpha}}\right)$ to the ion momentum equation, so the ions are accelerated down the electron-pressure gradient. This is exactly the divergence of the pressure flux $nT_e/m_{\alpha}$ appearing in the boundary term of the momentum balance \eqref{eq:ionfluidsystem}. And, its discretization ---the second term of equation \eqref{eq:local_ion_momentum_eqn} (to be discussed shortly), which equals $-[(nT_e)_{i+1}-(nT_e)_{i-1}]/(2\Delta x)$ once $q_e=-1$ is inserted--- is consistent with the centered discrete electric field used in the kinetic solve.
\end{rem}


\section{The numerical scheme: the kinetic solver}

When solving equation \eqref{eq:VFP1d2v}, the Vlasov term is evolved explicitly, and the acceleration and collision terms are evolved implicitly. Implicitly discretizing the collision operators requires the macroscopic quantities $n,\upar,\Ta,\Te$ at the future times, which we compute by solving the coupled fluid system \eqref{eq:fluidelectron} and \eqref{eq:ionfluidsystem}. We use implicit-explicit Runge-Kutta (IMEX RK) methods for the time-stepping, alternating between solving the fluid system and kinetic equation at each stage. This section describes the various components to the VFP integrator: discretizing equation \eqref{eq:VFP1d2v} in phase space, the rank-adaptive integrator and conservative truncation procedure extended to cylindrical coordinates, and the fluid solver.

\subsection{Discretizing the kinetic ion equation in phase space}\label{sec:kinetic_disc}

We discretize the physical domain $\Omega_x=[a,b]$ using $N_x+1$ uniformly distributed nodes
\begin{equation}
    a = x_{\frac{1}{2}}<x_{\frac{3}{2}}<...<x_{N_x-\frac{1}{2}}<x_{N_x+\frac{1}{2}}=b,
\end{equation}
where $x_i$ is the center of cell $I_i=[x_{i-\frac{1}{2}},x_{i+\frac{1}{2}}]$, and $\Delta x=|I_i|$ for all $i=1,...,N_x$. Similarly, we discretize the velocity domain $\Omega_v=[0,R_{\perp}]\times[-R_{||},R_{||}]$ using $N_{\perp}$ and $N_{||}$ uniform cells in each dimension, respectively. Defining the cell-centered grids
\begin{equation}
    \mathbf{v}_{\perp}: v_{\perp,1}<...<v_{\perp,{N_{\perp}}},\qquad \mathbf{v}_{||}: v_{||,1}<...<v_{||,{N_{||}}},
\end{equation}
we get the tensor product of uniform computational grids, $\mathbf{v}_{\perp}\otimes\mathbf{v}_{||}$. We choose $R_{\perp}$ and $R_{||}$ large enough for sufficient decay of the solution. The width of each cell in the perpendicular and parallel directions are respectively denoted by $\Delta_{\perp}$ and $\Delta_{||}$. Note that $v_{\perp,1}-\frac{\Delta_{\perp}}{2}=0$ and $v_{\perp,{N_{\perp}}}+\frac{\Delta_{\perp}}{2}=R_{\perp}$.

The scheme goes as follows: discretize equation \eqref{eq:VFP1d2v} in the physical domain, and solve a 2V system at each spatial node for $f_i(\vperp,\vpar,t)\approx f(x=x_i,\vperp,\vpar,t)$ at each $i=1,...,N_x$. In doing so, we can leverage rank-adaptive solvers in the velocity domain since solutions to the LBFP and VFP equations are known to be low-rank in $\mathbf{v}$ due to their Maxwellian-like distributions, especially near equilibrium. Moreover, this can be done in parallel at each spatial node $x_i$ to improve the overall runtime. Let $n_i,u_{||,i},T_{\alpha,i},T_{e,i}$ denote the nodal values of the macroscopic quantities. We use an upwind discretization for the flux term, and centered differencing for the electric field,
\begin{equation}\label{eq:VFP_semidiscrete_x}
    \frac{\partial f_i}{\partial t} + \max(\vpar,0)\frac{f_i-f_{i-1}}{\Delta x} + \min(\vpar,0)\frac{f_{i+1}-f_{i}}{\Delta x} + \frac{q_{\alpha}}{m_{\alpha}}E_{||,i}\frac{\partial f_i}{\partial \vpar} = C_{\alpha\alpha}(f_i) + C_{\alpha e}(f_i),
\end{equation}
where
\begin{equation}
    E_{||,i} = \frac{1}{q_en_i}\frac{(nT_e)_{i+1}-(nT_e)_{i-1}}{2\Delta x},\qquad i=1,...,N_x.
\end{equation}

Next, we consider the matrix solution of $f_i$ discretized over $\mathbf{v}_{\perp}\otimes\mathbf{v}_{||}$, denoted by $\mathbf{f}_i(t)\in\mathbb{R}^{N_{\perp}\times N_{||}}$. A second-order upwind finite difference method is used to discretize the acceleration term to get $\mathbf{f}_i\mathbf{D}_i^{T}$, where the boundary condition in $\vpar$ assumes the distribution function decays to zero as $|\vpar|\rightarrow\infty$.
\begin{equation}\label{eq:Emat}
    \mathbf{D}_i = \frac{1}{2\Delta_{||}}\left[\max\left(\frac{q_{\alpha}E_{||,i}}{m_{\alpha}},0\right)\text{lowband}(1,-4,3) + \min\left(\frac{q_{\alpha}E_{||,i}}{m_{\alpha}},0\right)\text{upband}(-3,4,-1)\right],
\end{equation}
where lowband(1,-4,3) is the banded matrix with entries (1,-4,3) on the second subdiagonal, first subdiagonal, and main diagonal, respectively. Similarly, upband(-3,4,-1) is the banded matrix with entries (-3,4,-1) on the main diagonal, first superdiagonal, and second superdiagonal, respectively.

Fokker-Planck operators are commonly discretized using the Chang-Cooper method \cite{ChangC1970practical}. The Chang-Cooper method is preferred among practitioners because it is easy to implement, second-order accurate, equilibrium-preserving, relative entropy dissipative, mass conservative, and positivity-preserving (under a loose restriction on $\Delta t$); it was recently extended to more general nonlinear Fokker-Planck equations in \cite{PareschiZ2018structure}. The Fokker-Planck operator \eqref{eq:LBFP} applied to $f_i$ can be expressed in azimuthally-symmetric cylindrical coordinates $(\vperp,\vpar)$,
\begin{equation}
    C_{\alpha\beta,i} = \nu_{\alpha\beta,i}\underbrace{\frac{1}{\vperp}\frac{\partial}{\partial\vperp}\left[\vperp\left(\frac{T_{\beta,i}}{m_{\alpha}}\frac{\partial f_i}{\partial\vperp} + (\vperp-u_{\beta,\perp,i})f_i\right)\right]}_{\eqqcolon C_{\alpha\beta,i}^{\perp}} + \nu_{\alpha\beta,i}\underbrace{\frac{\partial}{\partial\vpar}\left[\frac{T_{\beta,i}}{m_{\alpha}}\frac{\partial f_i}{\partial\vpar} + (\vpar-u_{\beta,||,i})f_i\right]}_{\eqqcolon C_{\alpha\beta,i}^{||}}.
\end{equation}

In our case, $u_{\beta,\perp,i}=0$ and $u_{\beta,||,i}=u_{||,i}$. The Chang-Cooper method uses a weighted sum of the cell-centered values of the solution to approximate the drift at the cell boundaries; it uses a centered difference to approximate the diffusive term. For brevity, we only specify which functions to use in the Chang-Cooper discretization, since the original paper clearly outlines the method. The original formulation considers operators of the form
\begin{equation}\label{eq:ChangCooper}
    \frac{1}{A(v)}\frac{\partial}{\partial v}\left[B(v,t)f+C(v,t)\frac{\partial f}{\partial v}\right].
\end{equation}

Discretizing $C_{\alpha\beta}^{\perp}$ has $A(\vperp)=\vperp$, $B(\vperp,t)=\vperp^2$, and $C(\vperp,t)=\frac{T_{\beta}}{m_{\alpha}}\vperp$. Similarly, $C_{\alpha\beta}^{||}$ has $A(\vpar)=1$, $B(\vpar,t)=\vpar-u_{||}$, and $C(\vpar,t)=\frac{T_{\beta}}{m_{\alpha}}$. We refer the reader to \cite{ChangC1970practical} for more details, or to \cite{KahzaCTQH2026structure} where Chang-Cooper type discretizations are presented for cylindrical coordinates. The resulting discretization of the collision operator $C_{\alpha\alpha}(f_i)+C_{\alpha e}(f_i)$ can be compactly written as $\mathbf{C}_i^{\perp}\mathbf{f}_i + \mathbf{f}_i\mathbf{C}_i^{||,T}$, where the differentiation matrices corresponding to Chang-Cooper method are denoted by
\begin{equation}\label{eq:Cmats}
    \mathbf{C}_i^{\perp} = \nu_{\alpha\alpha,i}\mathbf{C}_{\alpha\alpha,i}^{\perp} + \nu_{\alpha e,i}\mathbf{C}_{\alpha e,i}^{\perp},\qquad \mathbf{C}_i^{||} = \nu_{\alpha\alpha,i}\mathbf{C}_{\alpha\alpha,i}^{||} + \nu_{\alpha e,i}\mathbf{C}_{\alpha e,i}^{||}.
\end{equation}

Using equations \eqref{eq:Emat} and \eqref{eq:Cmats} to further discretize equation \eqref{eq:VFP_semidiscrete_x} in velocity, we get the semi-discrete equation
\begin{equation}\label{eq:semidiscrete}
    \frac{d\mathbf{f}_i}{dt} + \frac{1}{\Delta x}(\mathbf{f}_i-\mathbf{f}_{i-1})*\max(\mathbf{v}_{||},\mathbf{0})^T + \frac{1}{\Delta x}(\mathbf{f}_{i+1}-\mathbf{f}_{i})*\min(\mathbf{v}_{||},\mathbf{0})^T + \mathbf{f}_i\mathbf{D}_i^T = \mathbf{C}_i^{\perp}\mathbf{f}_i + \mathbf{f}_i\mathbf{C}_i^{||,T},
\end{equation}
where $*$ denotes the Hadamard/element-wise product.

\begin{rem}
The Chang-Cooper discretization reduces to the second-order centered difference for the discrete Laplacian if $B=0$ and $A$ and $C$ are constants in \eqref{eq:ChangCooper}.
\end{rem}

\subsection{A rank-adaptive integrator in cylindrical coordinates}\label{sec:kinetic_rail}

Naively storing the entire tensor solution requires $\mathcal{O}(N_xN_{\perp}N_{||})$ degrees of freedom. But, the storage complexity can be significantly reduced if the solution exhibits low-rank structure. Just like how the singular value decomposition (SVD) can be used to give the best low-rank matrix approximation, there are several different matrix (tensor) decompositions that reduce the storage requirements while still giving an approximation within a small tolerance $\vartheta$. It is well known that solutions to Fokker-Planck, Vlasov-Fokker-Planck, and other kinetic equations have low-rank structure in velocity, especially near equilibrium. For our purposes, we assume the numerical solution $\mathbf{f}_i(t)$, with $1\leq i\leq N_x$, can be decomposed (up to a small tolerance $\vartheta$) as
\begin{equation}\label{eq:lowrank}
    \mathbf{f}_i(t) \approx \bfV_i^{\perp}(t)\mathbf{S}_i(t)\bfV_i^{||,T}(t),
\end{equation}
where $\bfV_i^{\perp}\in\mathbb{R}^{N_{\perp}\times r_i(t)}$, $\bfV_i^{||}\in\mathbb{R}^{N_{||}\times r_i(t)}$, $\mathbf{S}_i\in\mathbb{R}^{r_i(t)\times r_i(t)}$, and $r_i\ll\min\{N_{\perp},N_{||}\}$. Here, $\bfV_i^{\perp}$ and $\bfV_i^{||}$ are orthonormal one-dimensional bases in $\vperp$ and $\vpar$ respectively, $\mathbf{S}_i$ are the coefficients (also known as the core tensor/matrix), and $r_i(t)$ is the solution rank. The SVD is a special case of solution \eqref{eq:lowrank}, but in general $\mathbf{S}_i$ can be dense and contain negative values. This is a two-sided factorization (the matrix analogue of a Tucker decomposition), although there are other low-rank matrix/tensor decompositions that can be used. While $\mathbf{S}_i$ does not necessarily have to be square, we assume it is since we do not expect the rank in $\vperp$ and $\vpar$ to wildly vary. Notice that the storage complexity of the solution matrix has reduced to $\mathcal{O}(N_xN_vr)$ if $N_v=N_{\perp}=N_{||}$ and $r=r_1=...=r_{N_x}$ for the sake of analyzing the complexity. Even for three-dimensional simulations, such as the 1D-2V VFP equation, this can make a significant difference. The reduced computational cost of working with a low-rank solution, in addition to the parallelism that can be done in $x$, can significantly decrease the overall complexity of the scheme. Plugging the rank-adaptive solution \eqref{eq:lowrank} into equation \eqref{eq:semidiscrete},
\begin{equation}\label{eq:semidiscrete_lowrank}
    \frac{d}{dt}\Big(\bfV_i^{\perp}\mathbf{S}_i\bfV_i^{||,T}\Big) + \mathcal{F}_i + \bfV_i^{\perp}\mathbf{S}_i\Big(\mathbf{D}_i\bfV_i^{||}\Big)^T = \mathbf{C}_i^{\perp}\bfV_i^{\perp}\mathbf{S}_i\bfV_i^{||,T} + \bfV_i^{\perp}\mathbf{S}_i\Big(\mathbf{C}_i^{||}\bfV_i^{||}\Big)^T,
\end{equation}
where
\begin{align}
\begin{split}
    \mathcal{F}_i &= \frac{1}{\Delta x}\Big[\bfV_i^{\perp},\bfV_{i-1}^{\perp}\Big]
    \begin{bmatrix*}
        \mathbf{S}_i&0\\
        0&-\mathbf{S}_{i-1}
    \end{bmatrix*}
    \Big(\Big[\bfV_i^{||},\bfV_{i-1}^{||}\Big]*\max(\mathbf{v}_{||},\mathbf{0})\Big)^T\\
    &+ \frac{1}{\Delta x}\Big[\bfV_{i+1}^{\perp},\bfV_{i}^{\perp}\Big]
    \begin{bmatrix*}
        \mathbf{S}_{i+1}&0\\
        0&-\mathbf{S}_{i}
    \end{bmatrix*}
    \Big(\Big[\bfV_{i+1}^{||},\bfV_{i}^{||}\Big]*\min(\mathbf{v}_{||},\mathbf{0})\Big)^T.
\end{split}
\end{align}

When working in this rank-adaptive framework, the factor matrices $\bfV_i^{\perp}$, $\bfV_i^{||}$ and $\mathbf{S}_i$ are evolved rather than the entire matrix solution $\mathbf{f}_i$. We want to use an implicit-explicit method to solve equation \eqref{eq:semidiscrete_lowrank}. However, high-order implicit methods are generally more difficult to implement than explicit methods when dealing with low-rank decomposed solutions since they require knowledge of the future bases. We use the recently proposed Reduced Augmentation Implicit Low-Rank (RAIL) method from \cite{NakaoQE2025reduced} as the 2V solver at each spatial node $x_i$; when extending to 3V, one could consider using the three-dimensional extension of RAIL in \cite{NakaoCE2026low}. It is a rank-adaptive implicit(-explicit) integrator that naturally fits the framework of our proposed method and achieves second-order temporal accuracy or better. Our proposed method is meant to provide a general framework to solve the hybrid VFP equation, and any preferred high-order implicit-explicit rank-adaptive 2V solver can be used instead. We chose RAIL for its simplicity, but other similar high-order implicit(-explicit) rank-adaptive methods include \cite{MengAC2025preconditioning,KahzaTQC2024krylov,CerutiEKL2024robust,LiJ2026high}.

We solve the kinetic-ion equation using a second-order two-stage IMEX RK method since higher-order extensions are straightforward. The RAIL method can accommodate higher-order IMEX RK methods, with some Butcher tableaus found in \cite{AscherRS1997implicit}. The RAIL method shown in \cite{NakaoQE2025reduced} was presented in Cartesian coordinates. The extension to cylindrical coordinates requires a slight modification in order to compute the discrete weighted $\ell^2$ inner product.

Let $\mathbf{f}_i^k\approx \mathbf{f}_i(t^k)$ be the solution at time $t^k$ for $k=0,1,...,N_t$. To solve equation \eqref{eq:semidiscrete_lowrank}, we consider the following stiffly accurate second-order two-stage IMEX RK scheme (IMEX222):
\begin{subequations}
\begin{align}
\begin{split}\label{eq:IMEX222_stage1}
    \bfV_i^{\perp,(1)}\mathbf{S}_i^{(1)}\bfV_i^{||,(1),T} &= \bfV_i^{\perp,k}\mathbf{S}_i^k\bfV_i^{||,k,T} - \gamma\Delta t\mathcal{F}_i^{k} -\gamma\Delta t\bfV_i^{\perp,(1)}\mathbf{S}_i^{(1)}\Big(\mathbf{D}_i^{(1)}\bfV_i^{||,(1)}\Big)^T\\
    &+ \gamma\Delta t\mathbf{C}_i^{\perp,(1)}\bfV_i^{\perp,(1)}\mathbf{S}_i^{(1)}\bfV_i^{||,(1),T} + \gamma\Delta t\bfV_i^{\perp,(1)}\mathbf{S}_i^{(1)}\Big(\mathbf{C}_i^{||,(1)}\bfV_i^{||,(1)}\Big)^T,
\end{split}\\
\begin{split}\label{eq:IMEX222_stage2}
    \bfV_i^{\perp,k+1}\mathbf{S}_i^{k+1}\bfV_i^{||,k+1,T} &= \bfV_i^{\perp,k}\mathbf{S}_i^k\bfV_i^{||,k,T} - \delta\Delta t\mathcal{F}_i^{k} - (1-\delta)\Delta t\mathcal{F}_i^{(1)}\\
    &-(1-\gamma)\Delta t\bfV_i^{\perp,(1)}\mathbf{S}_i^{(1)}\Big(\mathbf{D}_i^{(1)}\bfV_i^{||,(1)}\Big)^T -\gamma\Delta t\bfV_i^{\perp,k+1}\mathbf{S}_i^{k+1}\Big(\mathbf{D}_i^{k+1}\bfV_i^{||,k+1}\Big)^T\\
    & + (1-\gamma)\Delta t\mathbf{C}_i^{\perp,(1)}\bfV_i^{\perp,(1)}\mathbf{S}_i^{(1)}\bfV_i^{||,(1),T} + \gamma\Delta t\mathbf{C}_i^{\perp,k+1}\bfV_i^{\perp,k+1}\mathbf{S}_i^{k+1}\bfV_i^{||,k+1,T}\\
    &+ (1-\gamma)\Delta t\bfV_i^{\perp,(1)}\mathbf{S}_i^{(1)}\Big(\mathbf{C}_i^{||,(1)}\bfV_i^{||,(1)}\Big)^T + \gamma\Delta t\bfV_i^{\perp,k+1}\mathbf{S}_i^{k+1}\Big(\mathbf{C}_i^{||,k+1}\bfV_i^{||,k+1}\Big)^T
\end{split}
\end{align}
\end{subequations}
where $\gamma=1-\sqrt{2}/2$, $\delta=1-1/(2\gamma)$, and superscripts denote the solution at stages $t^{(0)}=t^k$, $t^{(1)}=t^k+\gamma\Delta t$, and $t^{(2)}=t^{k+1}$.

The RAIL method solves for the factor matrices at each stage using a projection-based procedure: (1) project out the $\vpar$ dependence so that the system is only a function of $\vperp$ and $t$, and update $\bfV_i^{\perp}$; (2) project out the $\vperp$ dependence so that the system is only a function of $\vpar$ and $t$, and update $\bfV_i^{||}$; (3) with both bases updated, project out all velocity dependence and solve for $\mathbf{S}_i$; (4) truncate the updated solution according to a tolerance to adapt the rank. The low-rank subspace that we project onto in each dimension includes the bases computed from the previous RK stages to ensure a richer subspace.

Before continuing onward, we describe the changes that are required to extend the RAIL method (and other similar projection-based methods) to cylindrical coordinates. Since $\mathbf{V}_i^{\perp}$ is the basis used in the $\vperp$ direction, we want to use the weighted inner product $\langle A,B\rangle_{\vperp}\coloneqq \int_0^{\infty}{AB\vperp d\vperp}$. At the discrete level, the orthogonal projection of $\mathbf{A}$ onto $\bfV_i^{\perp}$ is then $\mathbf{V}_i^{\perp}(\mathbf{v}_{\perp}*\mathbf{V}_i^{\perp})^T\mathbf{A}$. Furthermore, the SVD and (reduced) QR factorization are two methods to compress the matrix, and generate orthonormal matrices for the column and/or row spaces. However, numerically these are done in the unweighted $\ell^2$ inner product. As such, we must multiply and divide by $\sqrt{\mathbf{v}_{\perp}}$ to enforce the weighted inner product in $\vperp$. We refer to these \textit{weighted} factorizations as wSVD ($\mathbf{U}_w\mathbf{\Sigma V}^T$) and wQR ($\mathbf{Q}_w\mathbf{R}$) throughout this paper.

We now outline the rank-adaptive implicit-explicit integrator to solve equations \eqref{eq:IMEX222_stage1}-\eqref{eq:IMEX222_stage2}, and refer the reader to \cite{NakaoQE2025reduced} for finer details on the RAIL scheme.

\subsubsection{Stage 1}
We start by solving the first stage of the IMEX RK method given by equation \eqref{eq:IMEX222_stage1}.

\textbf{K Step.} Let $\mathbf{K}=\mathbf{V}_i^{\perp}\mathbf{S}_i\mathbf{V}_i^{||,T}\mathbf{V}_i^{||,k}\in\mathbb{R}^{N_{\perp}\times r_i^k}$ be the projected solution onto $\bfV_i^{||,k}$. Projecting equation \eqref{eq:IMEX222_stage1} and the solution in $\vpar$ onto the current basis $\bfV_i^{||,k}$, we get the Sylvester equation
\begin{equation}\label{eq:Keqn}
    \Big(\mathbf{I}-\gamma\Delta t\mathbf{C}_i^{\perp,(1)}\Big)\bfK^{(1)} + \bfK^{(1)}\Big(\gamma\Delta t\big(\mathbf{D}_i^{(1)} - \mathbf{C}_i^{||,(1)}\big)\bfV_i^{||,k}\Big)^T\bfV_i^{||,k} = \bfK^k - \gamma\Delta t\mathcal{F}_i^k\bfV_i^{||,k}.
\end{equation}

After solving equation \eqref{eq:Keqn} (see \cite{Simoncini2016computational} for a comprehensive review of Sylvester solvers), an updated basis $\tilde{\tilde{\bfV}}_i^{\perp,(1)}$ can be extracted from $\bfK^{(1)}$ by computing a reduced QR factorization (using \texttt{qrw}), referring to the definition of $\bfK$. That is, $\bfK^{(1)}=\mathbf{Q}_w\mathbf{R}\eqqcolon\tilde{\tilde{\bfV}}_i^{\perp,(1)}\mathbf{R}$, where the double tilde is used to denote the update basis obtained from the K step.

\textbf{L Step.} Let $\mathbf{L}=\bfV_i^{||}\mathbf{S}_i^T\bfV_i^{\perp,T}\big(\mathbf{v}_{\perp}*\bfV_i^{\perp,k}\big)$ be the projected solution onto $\bfV_i^{\perp,k}$, where we have transposed the solution so that the column space is for $\vpar$. Repeating the same procedure as the K step, we project onto $\bfV_i^{\perp,k}$ to get the Sylvester equation
\begin{equation}\label{eq:Leqn}
    \Big(\mathbf{I}+\gamma\Delta t\big(\mathbf{D}_i^{(1)}-\mathbf{C}_i^{||,(1)}\big)\Big)\mathbf{L}^{(1)} + \mathbf{L}^{(1)}\Big(-\gamma\Delta t\mathbf{C}_i^{\perp,(1)}\bfV_i^{\perp,k}\Big)^T\big(\mathbf{v}_{\perp}*\bfV_i^{\perp,k}\big) = \mathbf{L}^k - \gamma\Delta t\mathcal{F}_i^{k,T}\big(\mathbf{v}_{\perp}*\bfV_i^{\perp,k}\big).
\end{equation}

After solving equation \eqref{eq:Leqn}, let $\mathbf{L}^{(1)}=\mathbf{QR}\eqqcolon\tilde{\tilde{\bfV}}_i^{||,(1)}\mathbf{R}$.

\textbf{S Step.} We perform a Galerkin projection to get a matrix equation for the core tensor $\mathbf{S}_i$. To ensure that the projection subspace possesses information from both the updated and previous bases, we augment the two together:
\begin{equation}\label{eq:augmatrix}
    \Big[\tilde{\tilde{\bfV}}_i^{\perp,(1)},\bfV_i^{\perp,k}\Big],\qquad \Big[\tilde{\tilde{\bfV}}_i^{||,(1)},\bfV_i^{||,k}\Big].
\end{equation}

Since the rank of the augmented bases has nearly doubled, the \textit{reduced augmentation procedure} is applied to reduce the rank, as per Algorithm C.1 in \cite{NakaoQE2025reduced}. In short, a reduced QR (or wQR) factorization is applied to each augmented matrix in \eqref{eq:augmatrix} to orthonormalize the columns; then the upper triangular matrices $\mathbf{R}$ are compressed using the truncated SVD (or wSVD) according to a tolerance of $10^{-14}$ on the singular values; then the $\mathbf{Q}$ (or $\mathbf{Q}_w$) is multiplied with the left singular vectors to get the updated bases $\hat{\bfV}_i^{\perp,(1)}$ and $\hat{\bfV}_i^{||,(1)}$. This tolerance of $10^{-14}$ is chosen so that it is small enough to not affect the consistency of the scheme, but large enough to remove redundant basis vectors.

Let $\hat{\mathbf{S}}_i=\big(\mathbf{v}_{\perp}*\hat{\bfV}_i^{\perp,(1)}\big)^T\bfV_i^{\perp}\mathbf{S}_i\bfV_i^{||,T}\hat{\bfV}_i^{||,(1)}$. Projecting equation \eqref{eq:IMEX222_stage1} and the solution onto the updated bases, we solve the Sylvester equation
\begin{align}\label{eq:Seqn}
\begin{split}
    \Big(\mathbf{I}-\gamma\Delta t\big(\mathbf{v}_{\perp}*\hat{\bfV}_i^{\perp,(1)}\big)^T\big(\mathbf{C}_i^{\perp,(1)}\hat{\bfV}_i^{\perp,(1)}\big)\Big)\hat{\mathbf{S}}_i^{(1)} + \hat{\mathbf{S}}_i^{(1)}&\Big(\gamma\Delta t\big(\big(\mathbf{D}_i^{(1)} - \mathbf{C}_i^{||,(1)}\big)\hat{\bfV}_i^{||,(1)}\big)^T\Big)\hat{\bfV}_i^{||,(1)}\\
    &= \hat{\mathbf{S}}_i^k - \gamma\Delta t\big(\mathbf{v}_{\perp}*\hat{\bfV}_i^{\perp,(1)}\big)^T\mathcal{F}_i^k\hat{\bfV}_i^{||,(1)}.
\end{split}
\end{align}

\textbf{T Step.} To ensure that the updated solution $\hat{\bfV}_i^{\perp,(1)}\hat{\mathbf{S}}_i^{(1)}\hat{\bfV}_i^{||,(1),T}$ remains as low-rank as possible, we perform one final truncation. Naively using the wSVD to compress the solution destroys conservation up to the tolerance $\vartheta$ that is used to truncate the singular values. To truncate the solution (according to tolerance $\vartheta$) while conserving the mass, momentum, and energy of the system, we use the Local Macroscopic Conservative (LoMaC) method in \cite{GuoQ2024local}, which can be viewed as a conservative truncation procedure. The LoMaC method has been used in many rank-adaptive papers \cite{nelson2026local} where conservation of the moments was important, but presented in Cartesian coordinates; we show the extension to azimuthally symmetric cylindrical coordinates in subsection \ref{sec:kinetic_lomac}. After applying LoMaC, the final updated solution at Stage 1 is $\bfV_i^{\perp,(1)}\mathbf{S}_i^{(1)}\bfV_i^{||,(1),T}$.

\subsubsection{Stage 2}

We follow the same KLST steps from Stage 1 for equation \eqref{eq:IMEX222_stage2}, but with a slightly different projection subspace. We now include the bases computed at Stage 1 in the projection, as well as a first-order prediction at time $t^{k+1}$. Doing so utilizes information from all stages of the RK method to generate an approximation space as large as possible to observe the high-order accuracy obtained in the final stage of the RK method; see \cite{NakaoQE2025reduced,NakaoCE2026low} for more details.

\textbf{K and L Steps.} We follow the same procedure described in the K and L steps in Stage 1, but we instead project using the augmented bases
\begin{equation}\label{eq:augbases_stage2}
    \Big[\tilde{\bfV}_i^{\perp,k+1},\bfV_i^{\perp,(1)},\bfV_i^{\perp,k}\Big],\qquad \Big[\tilde{\bfV}_i^{||,k+1},\bfV_i^{||,(1)},\bfV_i^{||,k}\Big],
\end{equation}
where $\tilde{\bfV}_i^{\perp,k+1}$ and $\tilde{\bfV}_i^{||,k+1}$ are first-order predictions computed using the same first-order IMEX RK method from Stage 1 (but with an entire time-step $\Delta t$ instead of $\gamma\Delta t$). Since the rank has increased from including three bases in our augmented matrices in \eqref{eq:augbases_stage2}, we apply the same reduced augmentation procedure from the previous S step; the resulting bases are denoted $\check{\bfV}_i^{\perp,k+1}$ and $\check{\bfV}_i^{||,k+1}$. Projecting equation \eqref{eq:IMEX222_stage2} and the solution in $\vpar$ onto $\check{\bfV}_i^{||,k+1}$ yields the same Sylvester equation \eqref{eq:Keqn}, but the righthand side of the equation now includes the additional terms computed at $t^{(1)}$; similarly for the analogous L equation when projecting with $\check{\bfV}_i^{\perp,k+1}$.

\textbf{S and T Steps.} The S and T steps also follow the same procedure described in Stage 1. For the augmented bases used to get $\hat{\bfV}_i^{\perp,k+1}$ and $\hat{\bfV}_i^{||,k+1}$, apply the reduced augmentation procedure on \eqref{eq:augbases_stage2}, but replace the first-order predictions with the updated bases $\tilde{\tilde{\bfV}}_i^{\perp,k+1}$ and $\tilde{\tilde{\bfV}}_i^{||,k+1}$ computed from the Sylvester equations in the K and L steps. The resulting Sylvester equation for $\hat{\mathbf{S}}_i^{k+1}$ is the same as equation \eqref{eq:Seqn}, but the righthand side now includes the additional terms computed at $t^{(1)}$, as per equation \eqref{eq:IMEX222_stage2}. Applying the conservative LoMaC truncation procedure yields the final updated solution $\bfV_i^{\perp,k+1}\mathbf{S}_i^{k+1}\bfV_i^{||,k+1,T}$.

\subsection{LoMaC in azimuthally symmetric cylindrical coordinates}\label{sec:kinetic_lomac}

Naively truncating the solution using the SVD only ensures conservation of the moments up to the truncation tolerance $\vartheta$. A post-processing step that has been used in other rank-adaptive methods is the Local Macroscopic Conservative (LoMaC) truncation procedure \cite{GuoQ2024local}, which we briefly summarize and present for azimuthally symmetric cylindrical coordinates. The core idea is to: (1) scale the distribution function by a sufficiently decaying weight function, (2) project the scaled distribution function onto the subspace that conserves the desired moments, (3) truncate the part of the distribution function that lies in the orthogonal complement. As mentioned in \cite{NakaoCE2026low}, we found that the truncation re-introduces a small amount of artificial mass, momentum, and energy back into the system due to numerical roundoff error. We project this non-physical part out. Letting $\mathbf{f}^{\star}$ be the pre-truncated solution, the resulting truncated solution has the form
\begin{equation}\label{eq:lomac}
    \mathbf{f} = \mathbf{f}_1^M + (I-P)(\mathcal{T}_{\vartheta}(\mathbf{f}_2)),
\end{equation}
where $\mathbf{f}_1^M$ is the part that conserves the macroscopic quantities, and $\mathbf{f}_2=(I-P)(\mathbf{f}^{\star})$ is the part with zero moments, and $\mathcal{T}_{\vartheta}$ is some truncation operator like QR-SVD factorization with respect to tolerance $\vartheta$. For the LBFP and VFP equations, we want to conserve mass, momentum, and energy. If, e.g., only mass, is conserved, then the subspace just needs to be modified accordingly.

Let $\langle\cdot,\cdot\rangle_{\mathbf{w}}$ denote the weighted $\ell^2$ inner product with respect to weight function $\mathbf{w}$, chosen to be a Gaussian distribution in velocity space that decays fast enough for integrability. Note that this is still done in cylindrical coordinates, so the Jacobian determinant $\mathbf{v}_{\perp}$ is also included. Letting $\tilde{\mathbf{f}}\coloneqq \mathbf{f}^{\star}/\mathbf{w}$, orthogonally projecting (with respect to $\langle\cdot,\cdot\rangle_{\mathbf{w}}$) onto the subspace
\begin{equation}
    \mathcal{N} \coloneqq \text{span}\Big\{\mathbf{1}\ ,\ \mathbf{v}_{||}\ ,\ \mathbf{v}_{\perp}^2+\mathbf{v}_{||}^2\Big\}
\end{equation}
results in
\begin{equation}
    P_{\mathcal{N}}(\tilde{\mathbf{f}}) = \frac{n^{\star}}{\norm{\mathbf{1}}^2_\mathbf{w}}\mathbf{1} + \frac{(nu_{||})^{\star}}{\norm{ \mathbf{v}_{||}}^2_\mathbf{w}}\mathbf{v}_{||} + \frac{2(nU)^{\star} - cn^{\star}}{\norm{ \mathbf{v}_{\perp}^2+\mathbf{v}_{||}^2 - c\mathbf{1}}^2_\mathbf{w}} (\mathbf{v}_{\perp}^2+\mathbf{v}_{||}^2-c\mathbf{1}),
\end{equation}
where $c=\frac{\langle \mathbf{1}, \mathbf{v}_{\perp}^2+\mathbf{v}_{||}^2 \rangle_\mathbf{w}}{\norm{\mathbf{1}}^2_\mathbf{w}}$ is the shift that renders the energy basis vector $\mathbf{v}_{\perp}^2+\mathbf{v}_{||}^2-c\mathbf{1}$ orthogonal to $\mathbf{1}$ (so that the mass and energy projections decouple), and the stars $\star$ denote the moments of the pre-truncated solution. After projecting onto the subspace the conserves the moments, we rescale by $\mathbf{w}$ to get $P(\mathbf{f}^{\star})=\mathbf{w}*P_{\mathcal{N}}(\tilde{\mathbf{f}})$. Since the weight function can be separated into one-dimensional Gaussians in the perpendicular and parallel directions, $\mathbf{w}=\mathbf{w}_{\perp}\otimes\mathbf{w}_{||}$, the projected solution can be written in its low-rank form
\begin{equation}\label{eq:lomac_f1}
    P(\mathbf{f}^*) = 
    \left(\mathbf{w}_{\perp}*\begin{bmatrix*}
        \mathbf{1}_{\perp}&\mathbf{v}_{\perp}^2
    \end{bmatrix*}\right)
    \begin{bmatrix*}
        \frac{n^{\star}}{\norm{\mathbf{1}}^2_\mathbf{w}} - c\frac{2(nU)^{\star} - cn^{\star}}{\norm{ \mathbf{v}_{\perp}^2+\mathbf{v}_{||}^2 - c\mathbf{1}}^2_\mathbf{w}} & \frac{(nu_{||})^{\star}}{\norm{ \mathbf{v}_{||}}^2_\mathbf{w}} & \frac{2(nU)^{\star} - cn^{\star}}{\norm{ \mathbf{v}_{\perp}^2+\mathbf{v}_{||}^2 - c\mathbf{1}}^2_\mathbf{w}}\\
        \frac{2(nU)^{\star} - cn^{\star}}{\norm{ \mathbf{v}_{\perp}^2+\mathbf{v}_{||}^2 - c\mathbf{1}}^2_\mathbf{w}}&0&0
    \end{bmatrix*}
    \left(
    \mathbf{w}_{||}*\begin{bmatrix*}
        \mathbf{1}_{||}&\mathbf{v}_{||}&\mathbf{v}_{||}^2
    \end{bmatrix*}
    \right)^T.
\end{equation}

The factorized form of $\mathbf{f}_2=(I-P)(\mathbf{f}^{\star})$ is obtained by augmenting/concatenating the low-rank representations of $\mathbf{f}^{\star}$ and $P(\mathbf{f}^{\star})$, and then compressed/truncated according to tolerance $\vartheta$ using the QR-SVD approach proposed in \cite{GuoQ2022low}. However, since the radial bases are orthogonal with respect to the weighted inner product, we must use the \textit{weighted factorization} wQR; we naturally refer to this as the wQR-SVD approach.

If we want to enforce the local mass, momentum and energy to be $n_M$, $(nu_{||})_M$ and $(nU)_M$, then $\mathbf{f}_1^M$ takes the form of equation \eqref{eq:lomac_f1} but replaces the $\star$ moments with the $M$ moments. Lastly, we compute the final solution \eqref{eq:lomac} using the wQR-SVD procedure, but we compress using a tolerance of $10^{-14}$ since the reason for doing so is only to orthonormalize the final bases. The original LoMaC formulation does not do this, but it is imperative in projection-based integrators such as the RAIL method since the resulting bases must be orthonormal. 

\begin{rem}
The original LoMaC method in \cite{GuoQ2024local} scales and rescales the truncated part by $\sqrt{\mathbf{w}}$. Dividing by $\sqrt{\mathbf{w}}$ caused numerical instabilities in our experiments since it is nearly zero at the boundary of the domain. This was also observed and commented on in \cite{NakaoCE2026low}, and following their approach removing this scaling resolved the issue.
\end{rem}

\subsection{The fluid solver}

In order to solve the kinetic-ion equation at each spatial node using the conservative rank-adaptive solver described in subsections \ref{sec:kinetic_disc}-\ref{sec:kinetic_lomac}, the macroscopic quantities are required at each RK stage. For instance, $n_i^{(1)},u_{||,i}^{(1)},T_{\alpha,i}^{(1)},T_{e,i}^{(1)}$ are needed to construct $\mathbf{D}_i^{(1)}$, $\mathbf{C}_i^{\perp,(1)}$ and $\mathbf{C}_i^{||,(1)}$. Equation \eqref{eq:ionfluidsystem} models the total moments, but the local ion moment equations are similarly found without integrating over $\Omega_x$. Discretizing the local ion moment equations in $x$,
\begin{align}
\begin{split}\label{eq:local_ion_mass_eqn}
    \frac{dn_i}{dt} &= -\frac{1}{\Delta x}\Big[(\widehat{nu_{||}})_{i+\frac{1}{2}} - (\widehat{nu_{||}})_{i-\frac{1}{2}}\Big],
\end{split}\\
\begin{split}\label{eq:local_ion_momentum_eqn}
    \frac{d(nu_{||})_i}{dt} &= -\frac{1}{\Delta x}\Big[\hat{S}_{i+\frac{1}{2}}-\hat{S}_{i-\frac{1}{2}}\Big] + \frac{q_{\alpha}}{2\Delta xm_{\alpha}q_e}\Big[(nT_e)_{i+1}-(nT_e)_{i-1}\Big],
\end{split}\\
\begin{split}\label{eq:local_ion_energy_eqn}
    \frac{d(nU)_i}{dt} &= -\frac{1}{\Delta x}\Big[\hat{Q}_{i+\frac{1}{2}}-\hat{Q}_{i-\frac{1}{2}}\Big] + \frac{q_{\alpha}u_{||,i}}{2\Delta xm_{\alpha}q_e}\Big[(nT_e)_{i+1}-(nT_e)_{i-1}\Big]\\
    &\qquad\qquad\quad+ \frac{3\sqrt{2m_e}}{m_{\alpha}^2T_{e,i}^{3/2}}\Big[n_i^2T_{e,i}-\frac{m_{\alpha}}{3}\Big(2n_i(nU)_i-(nu_{||})_i^2\Big)\Big],
\end{split}
\end{align}
where the relationship \eqref{eq:UTrelationship} was used to express $T_{\alpha}$ in the collision operator of the ion energy equation. The upwind numerical fluxes are
\begin{equation}\label{eq:numfluxes}
    \begin{pmatrix*}
        \widehat{nu_{||}}\\
        \hat{S}\\
        \hat{Q}
    \end{pmatrix*}_{i+\frac{1}{2}} = 2\pi\int_{-\infty}^{\infty}\int_0^{\infty}{\vpar\hat{f}_{i+\frac{1}{2}}
    \begin{pmatrix*}
        1\\
        \vpar\\
        (\vperp^2+\vpar^2) /2
    \end{pmatrix*}
    \vperp d\vperp d\vpar},
\end{equation}
where
\begin{equation}
    \hat{f}_{i+\frac{1}{2}} = 
    \begin{cases}
        f_i,&\vpar>0,\\
        f_{i+1},&\vpar\leq 0.
    \end{cases}
\end{equation}

The numerical fluxes in equation \eqref{eq:numfluxes} are computed discretely using a midpoint Riemann sum. In the low-rank framework, this can be done one dimension at a time if the integrand is separable. That is, integrate over each low-rank basis separately before multiplying with the core matrix. This procedure is called the \textit{Low-Rank Double Integral (LRDI)} function in \cite{GalindoNPQT2025nodal}. Discretizing the electron energy equation \eqref{eq:fluidelectron} in $x$,
\begin{align}\label{eq:electronenergy_semidiscrete}
\begin{split}
    \frac{d(nT_e)_i}{dt} =& -\frac{5}{3\Delta x}\Big[\hat{u}_{||,i+\frac{1}{2}}(\widehat{nT_e})_{i+\frac{1}{2}} - \hat{u}_{||,i-\frac{1}{2}}(\widehat{nT_e})_{i-\frac{1}{2}}\Big] + \frac{(nu_{||})_i}{3n_i\Delta x}\Big[(nT_e)_{i+1}-(nT_e)_{i-1}\Big]\\
    &+ \frac{3.2}{3\sqrt{2m_e}\Delta x^2}\Big[(T_{e,i}^{5/2}+T_{e,i+1}^{5/2})(T_{e,i+1}-T_{e,i}) - (T_{e,i-1}^{5/2}+T_{e,i}^{5/2})(T_{e,i}-T_{e,i-1})\Big]\\
    &+ \frac{2\sqrt{2m_e}}{m_{\alpha}T_{e,i}^{3/2}}\Big[\frac{m_{\alpha}}{3}\Big(2n_i(nU)_i-(nu_{||})_i^2\Big)-n_i^2T_{e,i}\Big],
\end{split}
\end{align}
where
\begin{equation}
    \hat{u}_{||,i+\frac{1}{2}} = \frac{u_{||,i} + u_{||,i+1}}{2},
    \qquad
    (\widehat{nT_e})_{i+\frac{1}{2}} = 
    \begin{cases}
        n_iT_{e,i},&\hat{u}_{||,i+\frac{1}{2}}>0,\\
        n_{i+1}T_{e,i+1},&\hat{u}_{||,i+\frac{1}{2}}\leq 0.
    \end{cases}
\end{equation}

Equations \eqref{eq:local_ion_mass_eqn}-\eqref{eq:local_ion_energy_eqn} and \eqref{eq:electronenergy_semidiscrete} and discretized using the same second-order IMEX RK method used to update the ion distribution function. Doing so gives a nonlinear system of equations for $\Big[n_i^{(1)},(nu_{||})_i^{(1)},(nU)_i^{(1)},T_{e,i}^{(1)}\Big]$ at the first RK stage, and $\Big[n_i^{k+1},(nu_{||})_i^{k+1},(nU)_i^{k+1},T_{e,i}^{k+1}\Big]$ at the second RK stage, both solved using a Newton method. For each $i=1,...,N_x$, the local ion mass $n_i$ can be computed first since equation \eqref{eq:local_ion_mass_eqn} does not depend on the other macroscopic quantities.

Suppose we have already computed all the values of $n_i^{(1)}$ at the first stage of the RK method. Let $\mathbf{y}=\Big[(nu_{||})_1,...,(nu_{||})_{N_x},(nU)_1,...,(nU)_{N_x},T_{e,1},...,T_{e,N_x}\Big]^T$. Discretizing equations \eqref{eq:local_ion_momentum_eqn}, \eqref{eq:local_ion_energy_eqn}, \eqref{eq:electronenergy_semidiscrete} using the first stage of the IMEX RK method and rearranging the system of equations into the form $\mathbf{0}=\mathbf{R}\coloneqq\Big[\mathbf{R}_{nu_{||}};\mathbf{R}_{nU};\mathbf{R}_{T_e}\Big]$, Newton's method gives
\begin{equation}\label{eq:Newton}
    -\mathbf{R}_{\ell} = \mathbf{J}_{\ell}\delta\mathbf{y}_{\ell},
\end{equation}
where $\ell$ is the iteration counter, $\mathbf{J}_{\ell}$ is the Jacobian of $\mathbf{R}_{\ell}$ with respect to $\mathbf{y}_{\ell}$, and $\delta\mathbf{y}_{\ell}=\mathbf{y}_{\ell+1}-\mathbf{y}_{\ell}$. We note that the Jacobian has a convenient block form
\begin{equation}
    \mathbf{J}_{\ell} = 
    \begin{bmatrix*}
        \mathbf{J}_{nu_{||},nu_{||}}&\mathbf{J}_{nu_{||},nU}&\mathbf{J}_{nu_{||},T_e}\\
        \mathbf{J}_{nU,nu_{||}}&\mathbf{J}_{nU,nU}&\mathbf{J}_{nU,T_e}\\
        \mathbf{J}_{T_e,nu_{||}}&\mathbf{J}_{T_e,nU}&\mathbf{J}_{T_e,T_e}
    \end{bmatrix*}_{\ell}.
\end{equation}

The stopping criterion for the Newton iteration is $\norm{\mathbf{R}_{\ell}}_{\infty}\leq\min\{5\times 10^{-12},(5\times 10^{-10})\norm{\mathbf{R}_0}_{\infty}\}$. It is also worth noting that since the values of $\mathbf{n}=[n_1,...,n_{N_x}]^T$ are already computed at time $t^{(1)}$, we replace $\mathbf{n}_{\ell}$ with $\mathbf{n}^{(1)}$ in equation \eqref{eq:Newton}. Solving for the macroscopic quantities at the second stage of the RK method follows similarly. The values of $\mathbf{n}^{k+1}$ are computed solving equation \eqref{eq:local_ion_mass_eqn} with the same discretization (from the second stage of the IMEX RK method) used in equation \eqref{eq:IMEX222_stage2}. Then, equations \eqref{eq:local_ion_momentum_eqn}, \eqref{eq:local_ion_energy_eqn}, \eqref{eq:electronenergy_semidiscrete} are discretized using the same second stage equation, and solved with Newton's method in equation \eqref{eq:Newton}; values of $\mathbf{n}_{\ell}$ are replaced with those of $\mathbf{n}^{k+1}$.

To ensure that the discrete moments of the collision operator vanish exactly, we apply the \textit{quadrature corrected moments (QCM)} procedure described in \cite{TaitanoCS2017equilibrium,GalindoNPQT2025nodal}. The moments computed by the fluid solver are slightly perturbed so that the discrete moments of Maxwellian distribution generated by the moments exactly matches the moments themselves. The nonlinear system of equations that enforces these constraints is solved using Newton's method, and we refer the reader to \cite{TaitanoCS2017equilibrium,GalindoNPQT2025nodal} for more details.

\subsection{The VFP integrator}

With all the components in place, we outline the kinetic-ion fluid-electron VFP solver. The moments are needed at the future times to discretize the collision operator and electric field, and the ion distribution function is needed to compute the kinetic fluxes in equation \eqref{eq:numfluxes}. Therefore, we need to alternate between solving the fluid equations and kinetic-ion equation. At each stage of the RK method, the fluid equations are updated; then, the updated macroscopic quantities are used to discretize the collision operator and electric field in order to update the kinetic-ion distribution function; then, the updated distribution function is used to compute the kinetic fluxes in the fluid equations for the subsequent stages. The scheme described in this paper is outlined in Algorithm \ref{algo:vfpsolver}, with extensions to other IMEX RK methods being straightforward.

\begin{algorithm}[t!]
    \caption{Rank-adaptive solver for the hybrid kinetic-ion fluid-electron VFP model}
    \label{algo:vfpsolver}
    {\bf Input:} $\hat{\bfV}^{\perp,k}_i$, $\hat{\mathbf{S}}_i^k$, $\hat{\bfV}_i^{||,k}$\\
    {\bf Output:} $\hat{\bfV}^{\perp,k+1}_i$, $\hat{\mathbf{S}}_i^{k+1}$, $\hat{\bfV}_i^{||,k+1}$
    \begin{algorithmic}[1]
        \State Compute the numerical fluxes \eqref{eq:numfluxes} at $t^k$ using a midpoint Riemann sum.
        \State Discretizing equations \eqref{eq:local_ion_mass_eqn}-\eqref{eq:local_ion_energy_eqn} and \eqref{eq:electronenergy_semidiscrete} at the first stage of the IMEX RK method, solve for $\mathbf{n}^{(1)}$ and $\mathbf{y}^{(1)}$ using Newton's method \eqref{eq:Newton}.
        \State Correct the moments using the QCM procedure \cite{TaitanoCS2017equilibrium,GalindoNPQT2025nodal}; let $\tilde{\mathbf{n}}^{(1)}$ and $\tilde{\mathbf{y}}^{(1)}$ be the correct moments.
        \State Use $\tilde{\mathbf{n}}^{(1)}$ and $\tilde{\mathbf{y}}^{(1)}$ to discretize the collision operator \eqref{eq:Cmats} at $t^{(1)}$. Use $\mathbf{n}^{(1)}$ and $\mathbf{y}^{(1)}$ to discretize the acceleration term \eqref{eq:Emat} at $t^{(1)}$.
        \State Solve the KLS equations \eqref{eq:Keqn}-\eqref{eq:Seqn} for the solution $\hat{\bfV}_i^{\perp,(1)}\hat{\mathbf{S}}_i^{(1)}\hat{\bfV}_i^{||,(1),T}$.
        \State Truncate according to tolerance $\vartheta$ using the LoMaC procedure to get $\bfV_i^{\perp,(1)}\mathbf{S}_i^{(1)}\bfV_i^{||,(1),T}$. Use $\mathbf{n}^{(1)}\equiv\mathbf{n}^{(1)}_M$ and $\mathbf{y}^{(1)}\equiv\mathbf{y}^{(1)}_M$ for the macroscopic quantities in the LoMaC procedure to ensure conservation.
        \Statex\textit{Note: Steps 4-6 can be done in parallel at each spatial node $x_i$, for $i=1,...,N_x$.}
        \State Compute the numerical fluxes \eqref{eq:numfluxes} at $t^{(1)}$ using a midpoint Riemann sum. 
        \State Repeat steps 2-6 at the second stage of the IMEX RK method to get $\bfV_i^{\perp,k+1}\mathbf{S}_i^{k+1}\bfV_i^{||,k+1,T}$.
    \end{algorithmic}
\end{algorithm}


\section{Numerical experiments}

In this section, we present results for the rank-adaptive scheme (using RAIL for the rank-adaptive integrator) on a suite of problems. The order of accuracy and conservation of the moments is numerically verified on the heat equation, LBFP equation, and VFP equation; rank plots are also included. In particular, we include results for a standing shock problem of the hybrid VFP model, as well as a weak Landau damping test for the hybrid Vlasov equation. As mentioned in \cite{NakaoQE2025reduced}, the first $r^0\ll N_v$ singular vectors/values of the initial condition need to span a large enough space for the first time-step; we set $r^0=10$. While the focus of this paper was on the hybrid VFP model, the diversity of numerical tests is intended to demonstrate applicability across heat, Fokker-Planck, and (collisional and collisionless) hybrid Vlasov regimes in azimuthally symmetric cylindrical coordinates.

\subsection{Heat Equation}
\begin{figure}[h!]
\begin{minipage}[b]{0.47\linewidth}
    \centering
    \includegraphics[width=\textwidth]{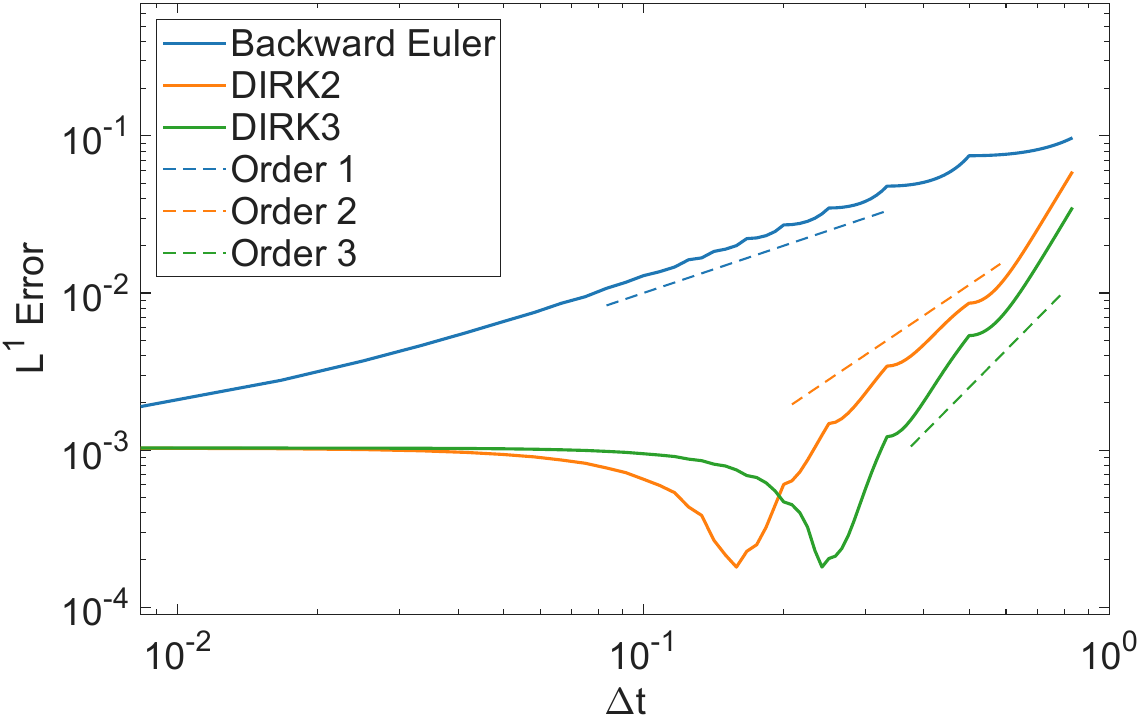}
    \caption{Error plot for equation \eqref{eq:heat-eqn} with solution \eqref{eq:heat-eqn-soln}.}
    \label{fig:heat-eqn-temporal-error-plot}
\end{minipage}
\begin{minipage}[b]{0.04\linewidth}
    \centering
    \ \\
\end{minipage}
\begin{minipage}[b]{0.47\linewidth}
    \centering
    \includegraphics[width=\textwidth]{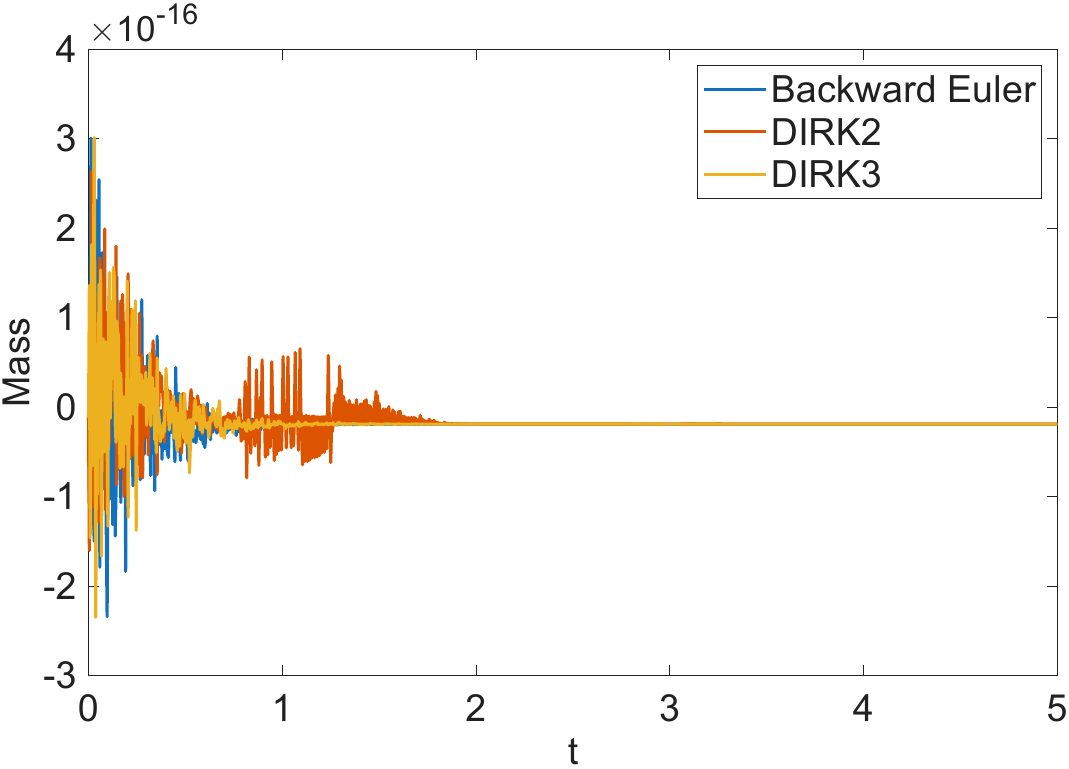}
    \caption{Mass plot for equation \eqref{eq:heat-eqn} with solution \eqref{eq:heat-eqn-soln}.}
    \label{fig:heat-eqn-mass}
\end{minipage}
\end{figure}

\begin{equation}
    \frac{\partial f}{\partial t}=\alpha \left(\frac{1}{\vperp}\frac{\partial}{\partial\vperp}\left(\vperp\frac{\partial f}{\partial\vperp}\right) + \frac{\partial^2 f}{\partial \vpar^2}\right), \qquad (\vperp,\vpar)\in(0,1)\times(0,1)
    \label{eq:heat-eqn}
\end{equation}
where as assume a homogeneous Neumann boundary condition at $\vperp=0$, homogeneous Dirichlet boundary condition at $\vperp=1$, and periodic boundary conditions in $\vpar$. Letting $\alpha=0.1$, we consider the first Fourier mode, with exact solution
\begin{equation}\label{eq:heat-eqn-soln}
    f(\vperp, \vpar, t)=\text{exp}\Big(-\alpha\big((j_{01})^2 +(2\pi)^2\big)t\Big)J_0(j_{01}\vperp)\sin(2\pi \vpar),
\end{equation}
where $J_0$ is the zeroth-order Bessel function of the first kind and $j_{01}$ is the first root of $J_0$. Equation \eqref{eq:heat-eqn} is solved using high-order stiffly accurate diagonally implicit Runge-Kutta (DIRK) methods \cite{AscherRS1997implicit}. The method follows identically to the IMEX RK presentation in subsection \ref{sec:kinetic_rail}, but with zero advection; extensions to higher-order RK methods is straightforward and described in \cite{NakaoQE2025reduced}. Only the kinetic solver is needed for this problem since the Laplacian does not require the updated moments of the solution. The purpose of this test is to demonstrate that the high-order temporal accuracy is observed while conserving mass and maintaining the rank-1 solution in cylindrical coordinates. As such, we use second-order centered differences for the spatial discretization. Figure \ref{fig:heat-eqn-temporal-error-plot} shows the expected accuracies in $L^1$ using a mesh of $N_{\perp}=40$, $N_{||}=80$, truncation tolerance $\vartheta=10^{-6}$, final time $T_f=1$, and LoMaC weight function $w=\text{exp}(-10(\vperp^2+\vpar^2))$. As expected, the temporal error of backward Euler dominates the spatial error even for small $\Delta t$. Whereas, there is an optimal time-stepping size for the high-order DIRK methods due to cancellation between the spatial and temporal errors. Noting that the mass of solution \eqref{eq:heat-eqn-soln} is zero, Figure \ref{fig:heat-eqn-mass} shows the mass is conserved up to machine precision. We report that the solution stayed rank-1 (occasionally jumping to rank-2), as it should based on \eqref{eq:heat-eqn-soln}. We used a mesh of $N_{\perp}=80$, $N_{||}=160$, time-step $\Delta t=1/(1/\Delta_{\perp}+1/\Delta_{||})$, final time $T_f=5$, and the same truncation tolerance and LoMaC weight function as before.

\subsection{Lenard-Bernstein-Fokker-Planck Equation}
Consider the spatially homogeneous 0D-2V LBFP equation in azimuthally symmetric cylindrical coordinates modeling the total effect of small-angle Coulomb collisions between ions of the same species.
\begin{equation}\label{eq:LBFP_ionion}
    \frac{\partial f}{\partial t} = C_{\alpha \alpha} = \nu_{\alpha\alpha}\nabla_v \cdot \left(\frac{T_{\alpha}}{m_{\alpha}}\nabla_v f + (\mathbf{v} - \mathbf{u}_{\alpha})f\right), \quad (\vperp,\vpar)\in (0,14)\times(-16,16),
\end{equation}
where we assume zero flux boundary conditions. We only need the kinetic solver since the system is assumed spatially homogeneous, and hence the macroscopic quantities are assumed constant. Scaling the ion mass to unity, the equilibrium solution is the Maxwellian distribution function
\begin{equation}\label{eq:LBFP_equilibrium}
    f_{\infty}(\vperp, \vpar)=\frac{n_{\alpha}}{(2\pi T_{\alpha})^{3/2}}\exp\left(-\frac{(\vperp - u_{\alpha,\perp})^2 + (\vpar - u_{\alpha,||})^2}{2T_{\alpha}}\right).
\end{equation}
We test the relaxation of the system using a rank-2 initial condition \eqref{eq:LBFP_IC} comprising two Maxwellians whose parameters are included in Table \ref{tab:dfp-IC-values}. We set $u_{\perp}=0$ so that drift only occurs in the $\vpar$ direction. The Chang-Cooper method is used to discretize the LBFP collision operator.
\begin{equation}\label{eq:LBFP_IC}
    f(\vperp, \vpar, t=0)=f_{M_1}(\vperp, \vpar) + f_{M_2}(\vperp, \vpar).
\end{equation}
\begin{table}[h!]
    \centering
    \begin{tabular}{|c|c|c|c|c|}
        \hline
        & $n$ & $u_{\perp}$ & $u_{||}$ & $T$\\
        \hline
        $f_{M_1}$ & $2.0$ & $0.0$ & $-0.5$ & $2.0$\\
        \hline
        $f_{M_2}$ & $1.0$ & $0.0$ & $0.9$ & $1.0$\\
        \hline
    \end{tabular}
    \caption{Macroscopic quantities that define the two Maxwellians in initial condition \eqref{eq:LBFP_IC}.}
    \label{tab:dfp-IC-values}
\end{table}

\begin{figure}[h!]
    \centering
    \includegraphics[width=0.5\linewidth]{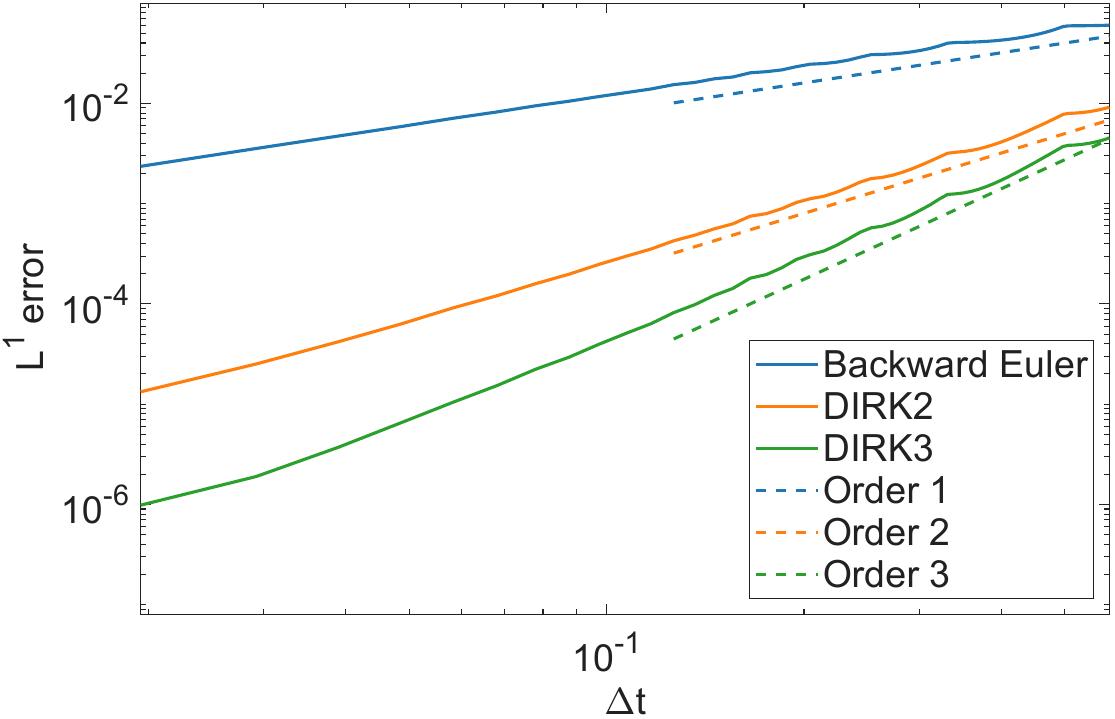}
    \caption{Error plot for equation \eqref{eq:LBFP_ionion} using backward Euler, DIRK2, and DIRK3. Mesh of $N_{\perp}=N_{||}=100$, final time $T_f=1$, and truncation tolerance $\vartheta=10^{-6}$.}
    \label{fig:dfp-temporal-error-plot}
\end{figure}

We numerically verify the temporal accuracy using a reference solution computed with mesh $N_{\perp}=N_{||}=100$, time-step $\Delta t=0.00974$, DIRK3, final time $T_f=1$, and truncation tolerance $\vartheta=10^{-6}$. As seen in Figure \ref{fig:dfp-temporal-error-plot}, varying the time-step shows the expected orders of accuracy for backward Euler, DIRK2, and DIRK3. We used a mesh of $N_{\perp}=N_{||}=100$, and tolerance $\vartheta=10^{-6}$.

\begin{figure}[h!]
\centering
    \begin{minipage}{0.49\linewidth}
        \centering
        \includegraphics[width=\linewidth]{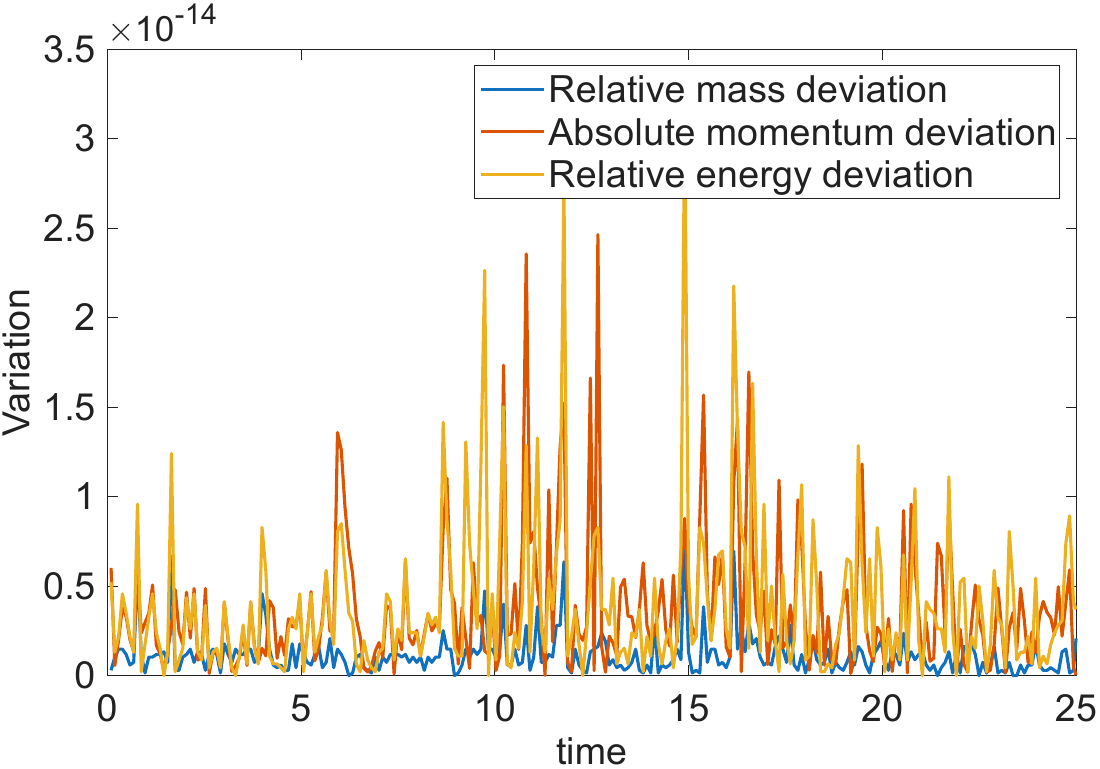} 
    \end{minipage}
    \hfill
    \begin{minipage}{0.49\linewidth}
        \centering
        \includegraphics[width=\linewidth]{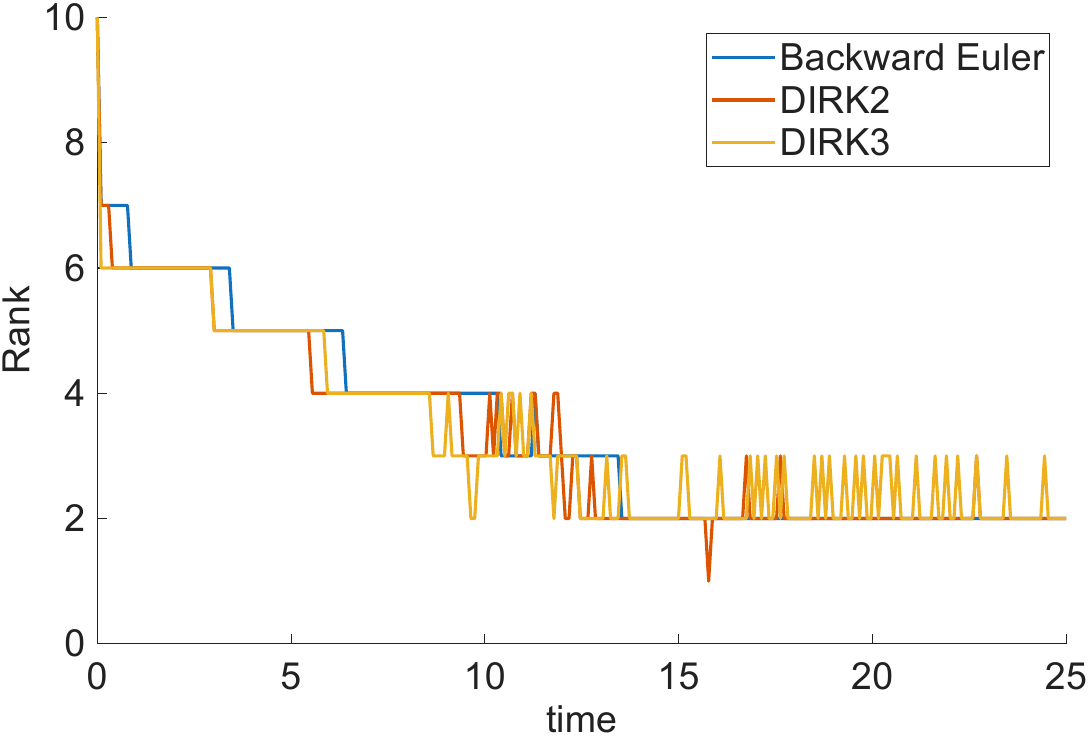}
    \end{minipage}
    \caption{Changes in the mass, momentum, and energy \textit{(left)}, and the solution rank \textit{(right)} for equation \eqref{eq:LBFP_ionion} with initial condition \eqref{eq:LBFP_IC}. Mesh $N_{\perp}=N_{||}=100$, truncation tolerance $\vartheta=10^{-6}$, time-step $\Delta t=1/(1/\Delta_{\perp}+1/\Delta_{||})$, and final time $T_f=25$. The changes in the mass, momentum, and energy are shown for DIRK3, but results were similar when using backward Euler and DIRK2.}
    \label{fig:dfp-conservation-rank-plot}
\end{figure}

\begin{figure}[h!]
\centering
    \begin{minipage}{0.49\linewidth}
        \centering
        \includegraphics[width=\linewidth]{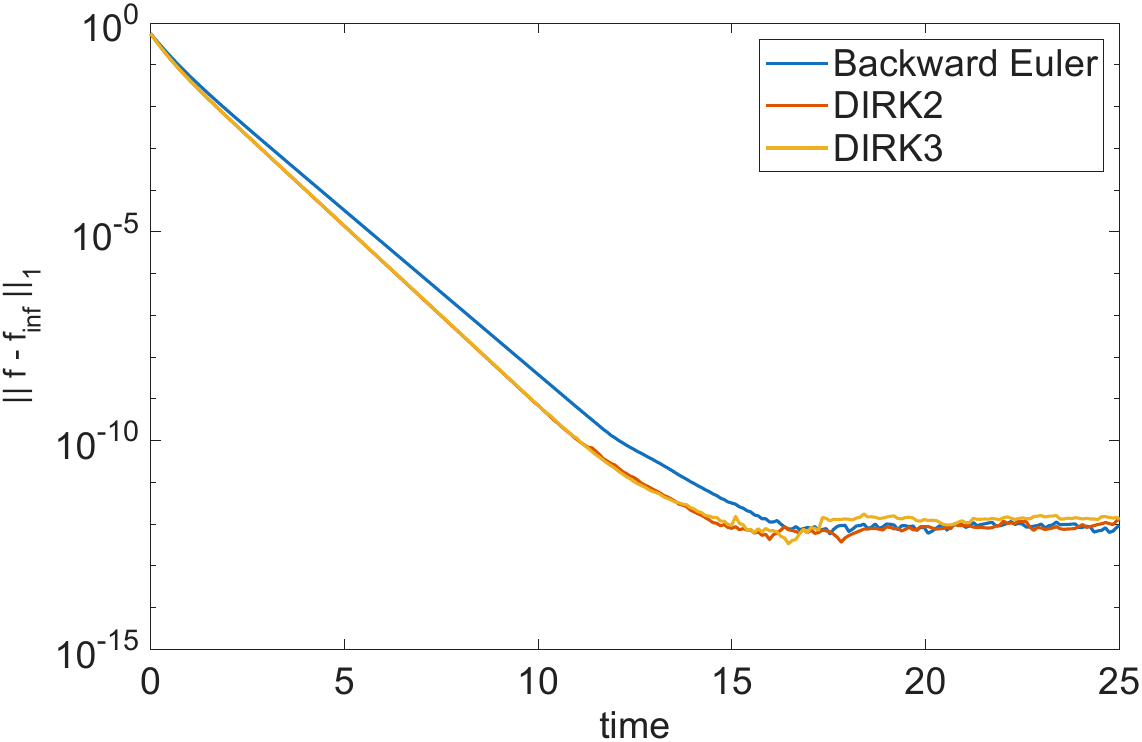}
    \end{minipage}
    \hfill
    \begin{minipage}{0.49\linewidth}
        \centering
        \includegraphics[width=\linewidth]{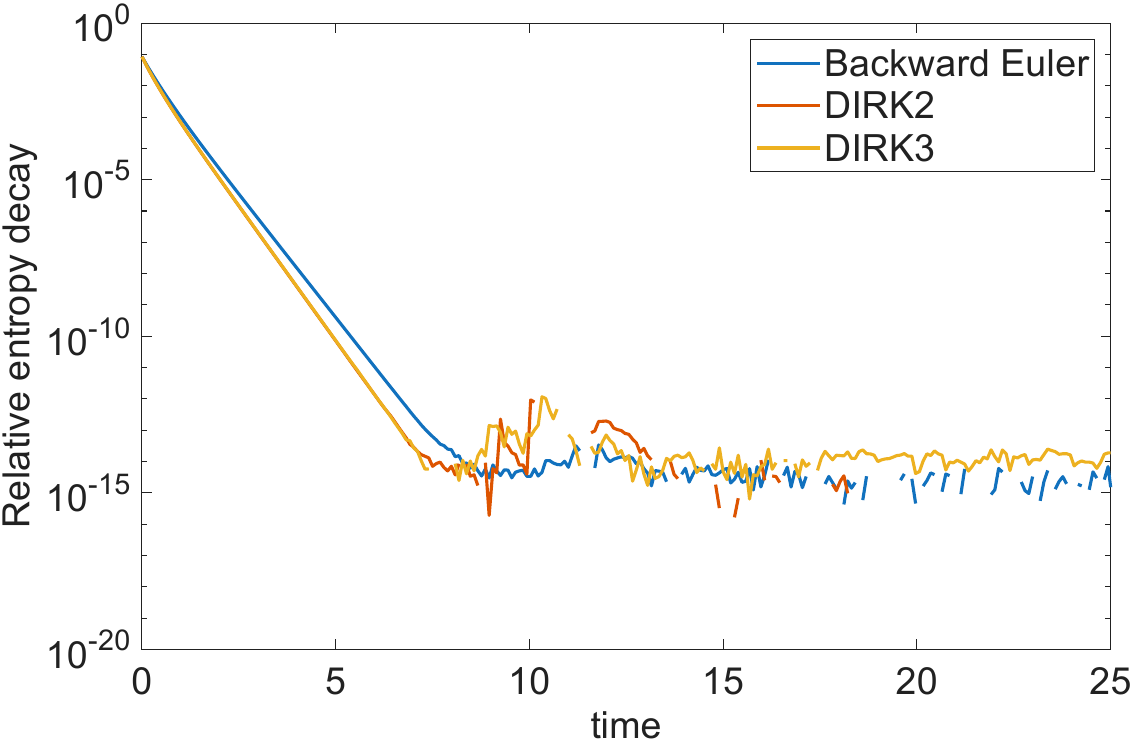}
    \end{minipage}
    \caption{Equilibrium preservation / $L^1$ decay towards relaxation \textit{(left)}, and Kullback relative entropy dissipation \textit{(right)} for equation \eqref{eq:LBFP_ionion} with initial condition \eqref{eq:LBFP_IC}. Mesh $N_{\perp}=N_{||}=100$, truncation tolerance $\vartheta=10^{-6}$, time-step $\Delta t=1/(1/\Delta_{\perp}+1/\Delta_{||})$, and final time $T_f=25$.}
    \label{fig:dfp-l1-relative-entropy}
\end{figure}

The solution rank and relative changes in the moments are shown in Figure \ref{fig:dfp-conservation-rank-plot}, and the $L^1$ decay towards relaxation $\norm{f-f_{\infty}}_1$ and Kullback relative entropy dissipation $\int_{\mathbb{R}^3}{f\log{(f/f_{\infty})} \vperp d^3v}$ of the solution are shown in Figure \ref{fig:dfp-l1-relative-entropy}. We use a mesh of $N_{\perp}=N_{||}=100$, truncation tolerance $\vartheta=10^{-6}$, time-step $\Delta t=1/(1/\Delta_{\perp}+1/\Delta_{||})$, and final time $T_f=25$. Figure \ref{fig:dfp-conservation-rank-plot} shows that the mass, momentum, and energy are conserved up to roughly machine precision. And, the rank-2 initial condition immediately jumps in rank as the two Maxwellians drift and diffuse before relaxing to the low-rank Maxwellian steady-state. Figure \ref{fig:dfp-l1-relative-entropy} also shows equilibrium preservation and relative entropy dissipation up to roughly machine precision, despite using a larger tolerance. This suggests that the LoMaC procedure projects onto a subspace that does a good job capturing the Maxwellian steady-state solution. It is imperative that the QCM procedure is used to compute the Maxwellian $f_{\infty}$ whose discrete moments match those of the initial condition \eqref{eq:LBFP_IC}. Otherwise, the equilibrium preservation and relative entropy dissipation only decay up to the numerical integration error which was roughly $\mathcal{O}(10^{-3})$.

\subsection{Standing shock problem}

We consider the 1D-2V hybrid Vlasov-Fokker-Planck model described by equations \eqref{eq:VFP1d2v}-\eqref{eq:fluidelectron}. We consider the Mach 5 standing shock problem tested in \cite{taitano_2018_jcp_r_adaptivity,TaitanoCS2021conservative,VidalMCL1993ion}, for which the steady-state solution is a standing shock. Since the shock front can continue to drift indefinitely when conservation is not enforced \cite{TaitanoCS2021conservative}, this serves as a good test for our proposed scheme. We show that our numerical scheme captures the standing shock behavior of the long-term solution while remaining mass, momentum, and energy conservative. Considering the domain $(x,\vperp,\vpar)\in(0,200)\times(0,8)\times(-8,8)$, the solution is initialized with a shock centered at $x=100$ defined by the downstream quantities $n_{\text{down}}=1$, $u_{||,\text{down}}=0.8676$, $T_{\alpha,\text{down}}=T_{e,\text{down}}=1$, and the upstream quantities $n_{\text{up}}=0.28$, $u_{||,\text{up}}=3.0984$, $T_{\alpha,\text{up}}=T_{e,\text{up}}=0.1152$. As in \cite{TaitanoCS2021conservative}, the initial shock is smoothed with hyperbolic tangents, with profiles given by
\begin{subequations}\label{eq:VFP_IC}
\begin{equation}
    n(x,t=0) = 0.36\tanh{(0.05(x-100))}+0.64,
\end{equation}
\begin{equation}
    u_{||}(x,t=0) = -1.1154\tanh{(0.05(x-100))} + 1.983,
\end{equation}
\begin{equation}
    T_{\alpha}(x,t=0) = 0.4424\tanh{(0.05(x-100))} + 0.5576.
\end{equation}
\end{subequations}

\begin{figure}[h!]
\centering
    \begin{minipage}{0.32\linewidth}
        \centering
        \includegraphics[width=\linewidth]{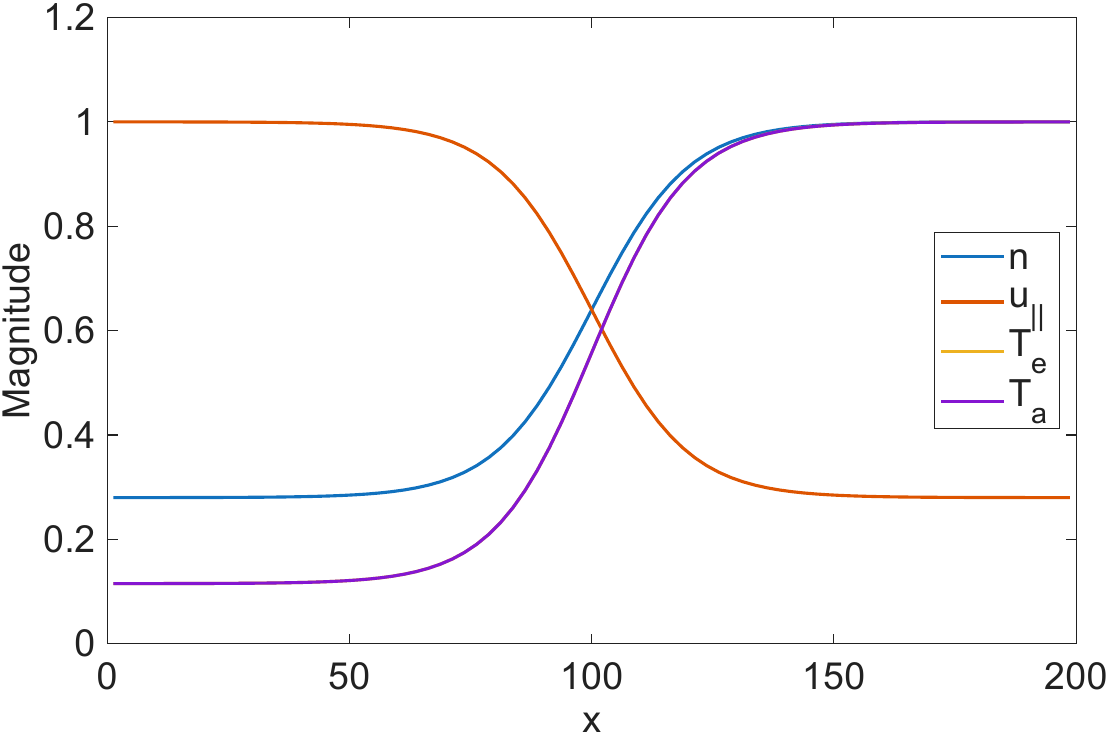}      
    \end{minipage}
    \hfill
    \begin{minipage}{0.32\linewidth}
        \centering
        \includegraphics[width=\linewidth]{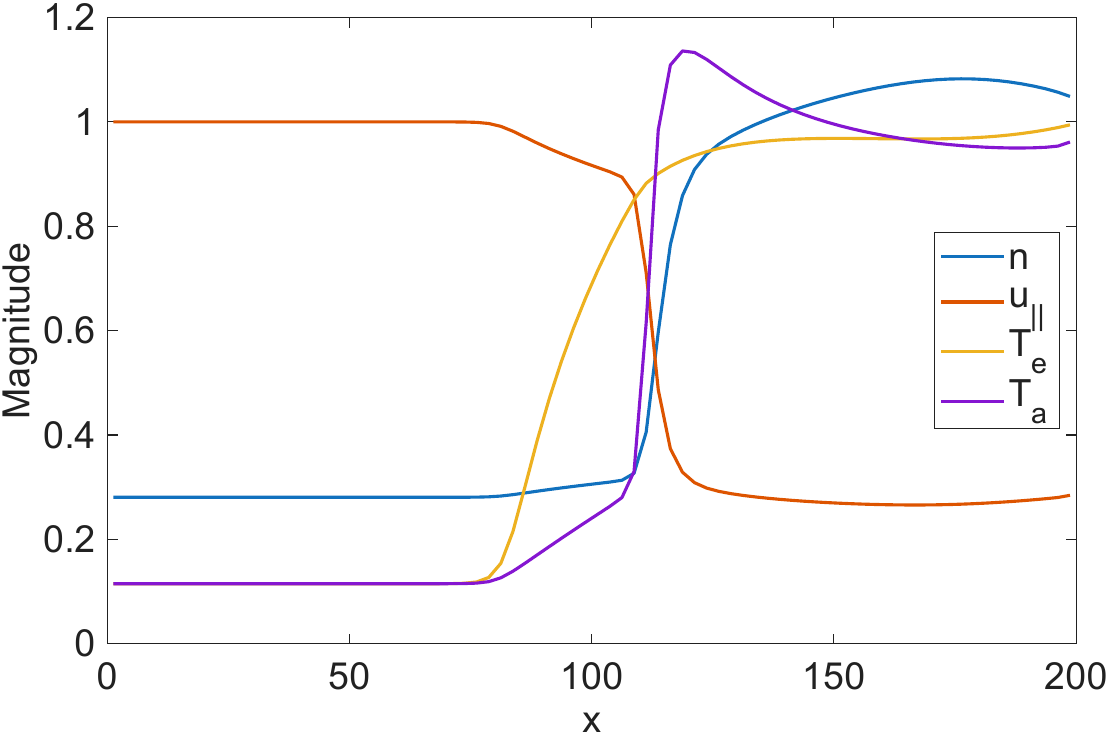}      
    \end{minipage}
    \hfill
    \begin{minipage}{0.32\linewidth}
        \centering
        \includegraphics[width=\linewidth]{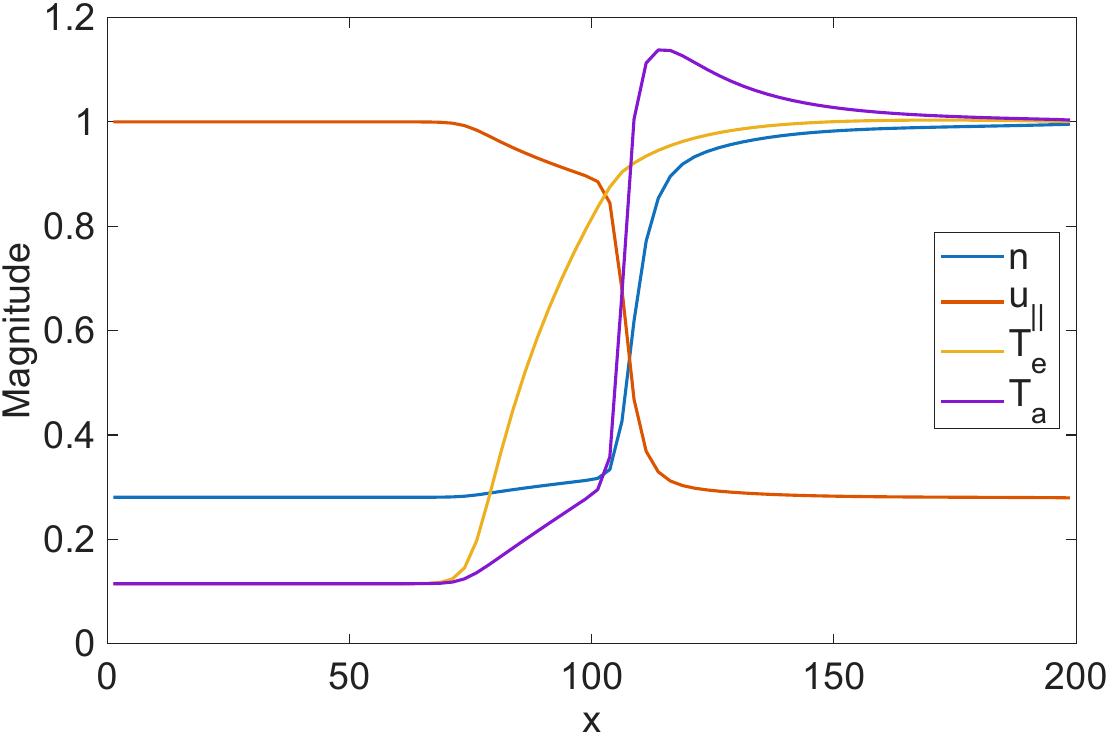}
    \end{minipage}
    \caption{Profiles of the macroscopic quantities at times $t=0$ \textit{(left)}, $t=100$, \textit{(middle)}, and $t=250$ \textit{(right)}. We used a mesh of $N_x=80$ and $N_{\perp}=N_{||}=100$, truncation tolerance $\vartheta=10^{-8}$, time-step $\Delta t=0.4$, and normalized $u_{||}$ to the upstream value.}
    \label{fig:vfp-macro-moments}
\end{figure}

The initial ion distribution function $f(x=x_i,\vperp,\vpar,t=0)$ is the Maxwellian distribution generated by the initial values of $n$, $u_{||}$ and $T_{\alpha}$ at the point $x_i$. We assume Dirichlet boundary conditions given by the downstream and upstream quantities. We set $\Delta t=0.005$ for the first step so that the initial dynamics can be captured, and use a fixed value of $\Delta t$ for all subsequent steps.

\begin{figure}[h!]
    \centering
    \includegraphics[width=0.6\linewidth]{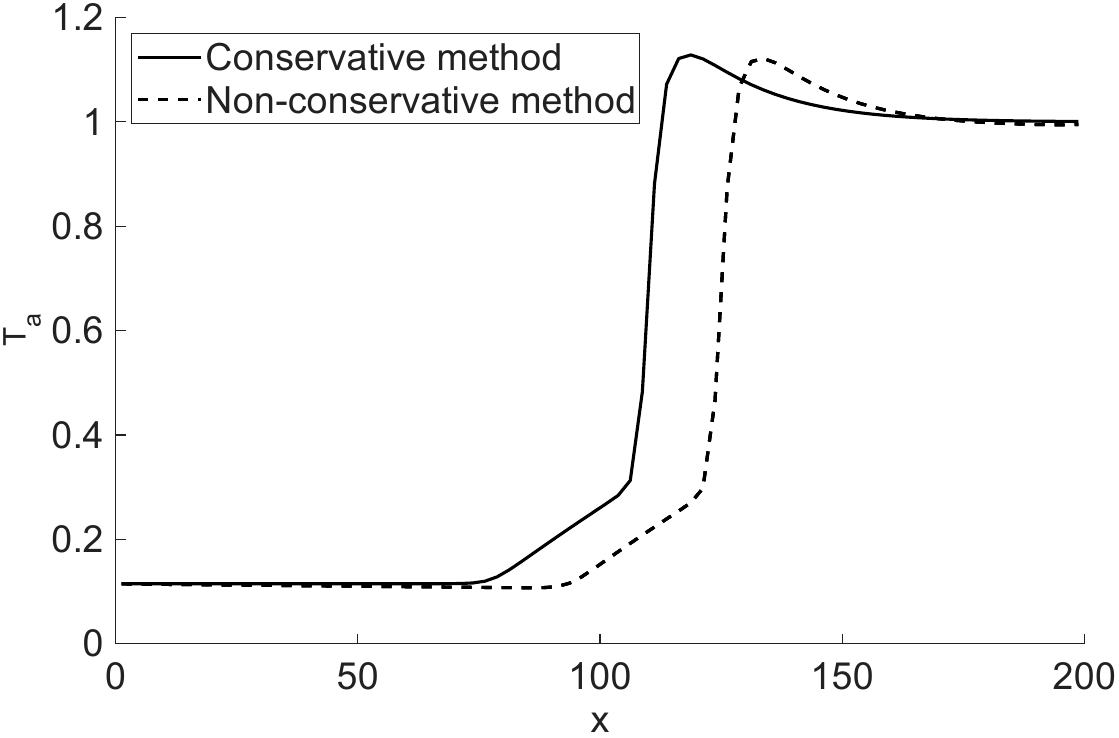}
    \caption{Ion temperature $T_\alpha$ profile with the conservative LoMaC truncation and QCM procedure for the Mach 5 standing shock problem \textit{(solid line)}. $T_\alpha$ profile simulated without using the LoMaC and QCM procedures \textit{(dashed line)}. Simulation run with IMEX222 on a mesh size of $N_x=80, N_{\perp}=N_{||}=100$, time-step $\Delta t = 0.4$ and truncation tolerance $\vartheta=10^{-8}$. $T_\alpha$ profiles are shown at final time $t=5000$. The shock drifts nonphysically when conservation is not enforced.}
    \label{fig:vfp_standing_shock_drift_test}
\end{figure}

\begin{figure}[h!]
\centering
    \begin{minipage}{0.49\linewidth}
        \centering
        \includegraphics[width=\linewidth]{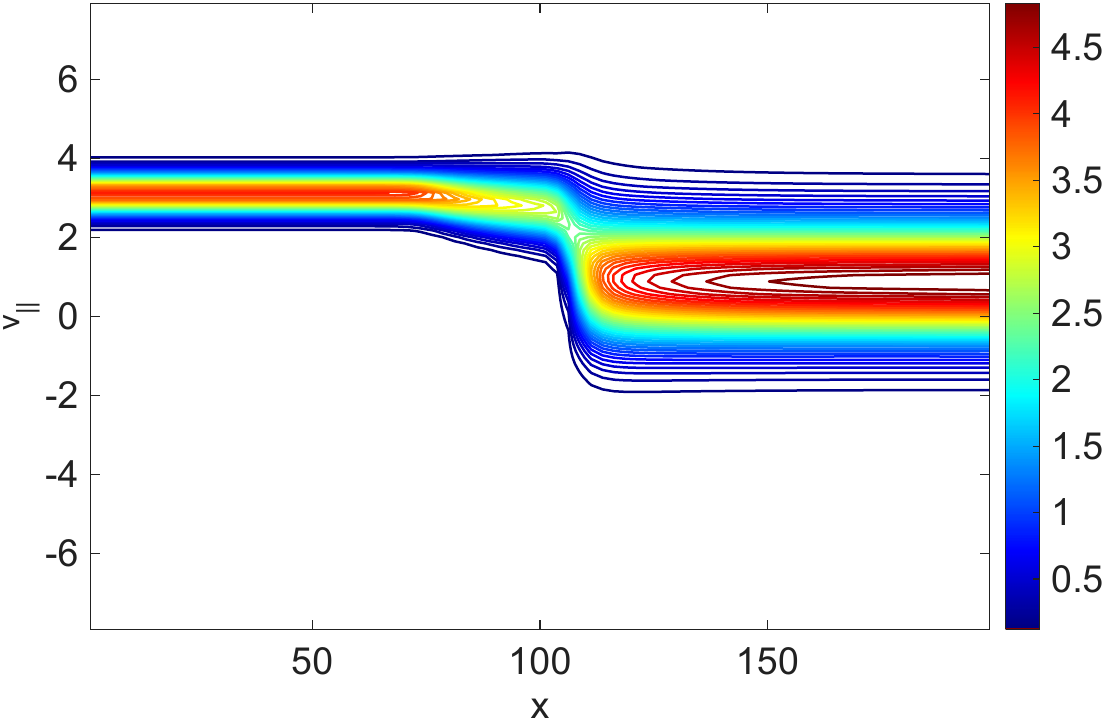}      
    \end{minipage}
    \hfill
    \begin{minipage}{0.49\linewidth}
        \centering
        \includegraphics[width=\linewidth]{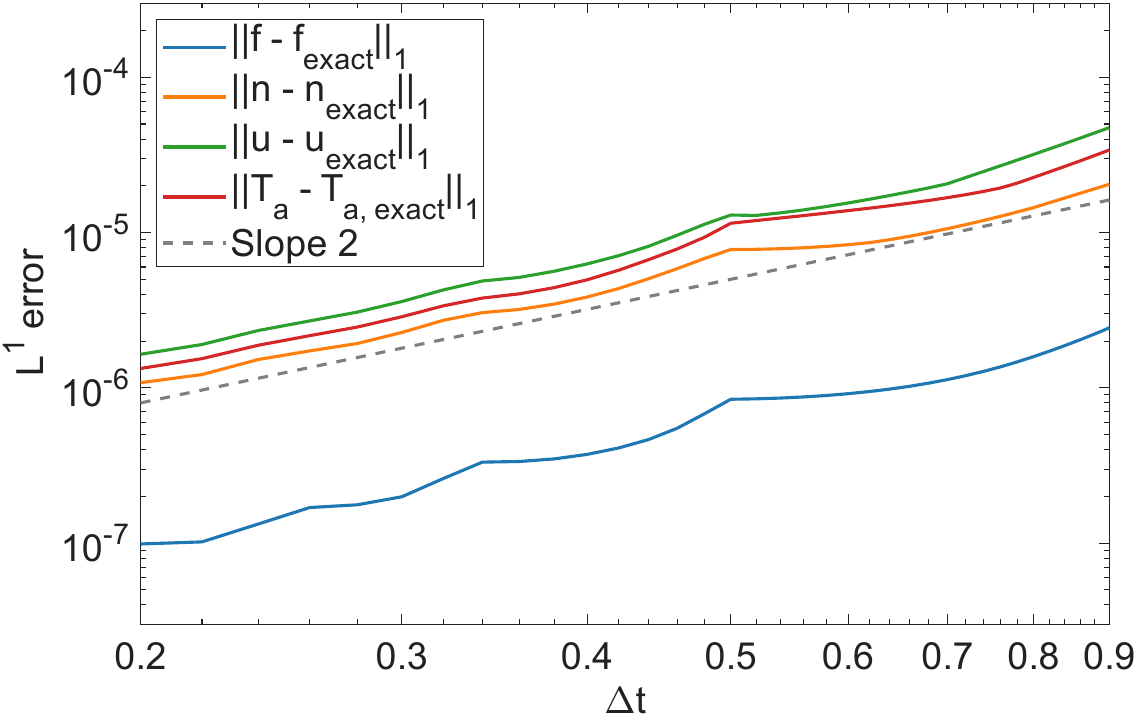}
    \end{minipage}
    \caption{The $\vperp$-integrated ion distribution function $2\pi\int_{\mathbb{R}_+}{f(x,\vperp,\vpar,t=250)\vperp d\vperp}$ for the Mach 5 steady-state shock, using mesh $N_x=80$, $N_{\perp}=N_{||}=100$, IMEX222, time-step $\Delta t=0.4$ \textit{(left)}. Error plot corresponding to IMEX222 with mesh $N_x=160$, $N_{\perp}=N_{||}=200$, final time $T_f=1$, and truncation tolerance $\vartheta=10^{-8}$ \textit{(right)}. The reference solution was computed with time-step $\Delta t=0.01$.}
    \label{fig:vfp-vperpint-temporal-error}
\end{figure}

\begin{figure}[h!]
\centering
    \begin{minipage}{0.49\linewidth}
        \centering
        \includegraphics[width=\linewidth]{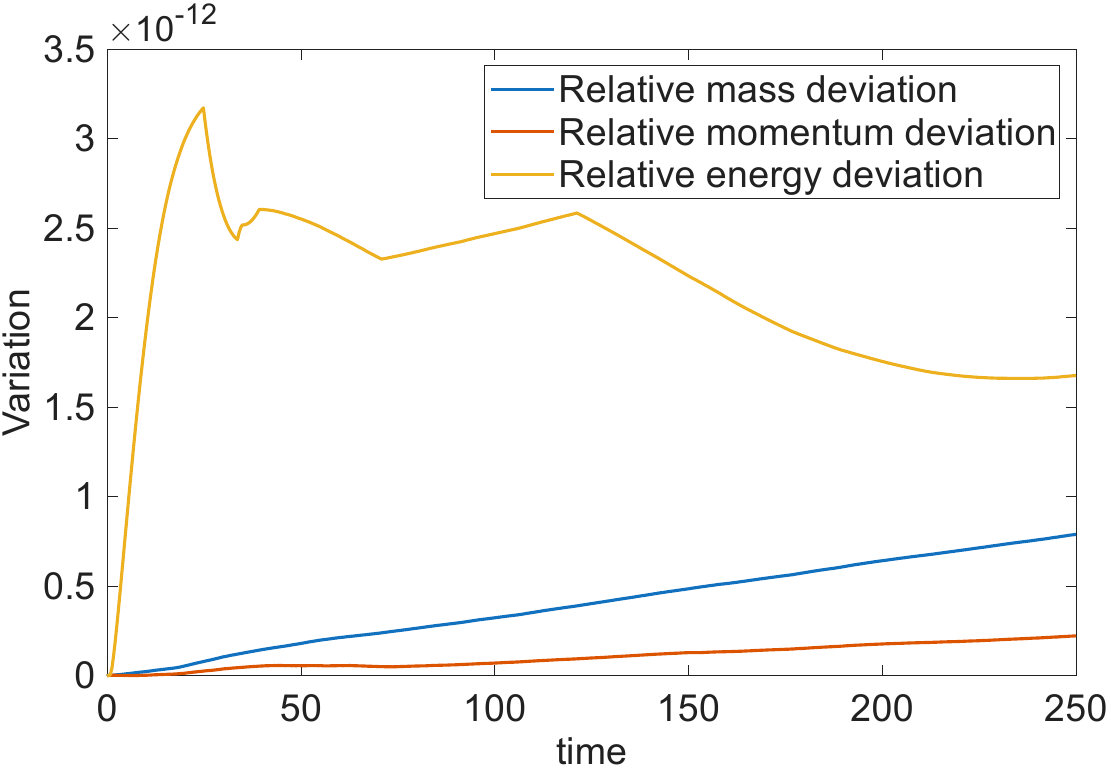}      
    \end{minipage}
    \hfill
    \begin{minipage}{0.49\linewidth}
        \centering
        \includegraphics[width=\linewidth]{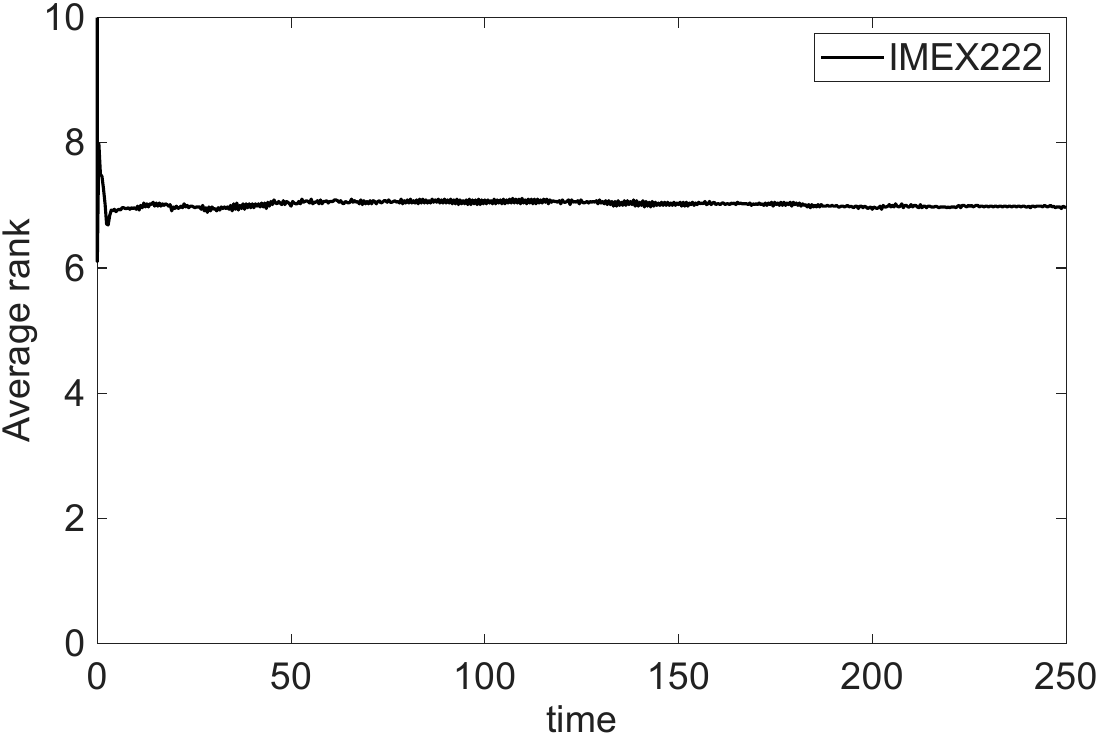}
    \end{minipage}
    \caption{Changes in mass, momentum, and energy \textit{(left)}, and the average solution rank across all spatial nodes $x_i$ \textit{(right)} for the Mach 5 steady-state shock. Mesh $N_x=80, N_{\perp}=N_{||}=100$, time-step $\Delta t = 0.4$, final time $T_f=250$, and truncation tolerance $\vartheta = 1.0E-8$. We note that when computing the relative change in the macroscopic quantities, we subtracted the contributions from the Dirichlet boundary conditions in equation \eqref{eq:ionfluidsystem}.}
    \label{fig:vfp-conservation-and-rank}
\end{figure}

Snapshots of the macroscopic quantities are shown in Figures \ref{fig:vfp-macro-moments} and \ref{fig:vfp_standing_shock_drift_test}. Notably, the ion temperature continues to drift to the right as $t$ increases when the QCM procedure and LoMaC truncation are not used. This demonstrates how critical conservation is for long-time simulations, without which we obtain nonphysical solutions. The $\vperp$-integrated distribution function is shown in Figure \ref{fig:vfp-vperpint-temporal-error}, in which we observe a clear elbow structure forming across the shock front at $x\sim 100$, similar to the results of Ref. \cite{taitano_2018_jcp_r_adaptivity}. However, unlike that work---which used the Rosenbluth-Fokker-Planck operator that supports a particle-velocity-dependent collision frequency---no energetic beam features form here. This difference is expected from the Lenard-Bernstein Fokker-Planck operator, since all particles effectively collide at the same frequency, relaxing the high-energy particles toward the Maxwellian faster than the more complete model predicts. Furthermore, the relative mass, momentum, and energy are shown in Figure \ref{fig:vfp-conservation-and-rank}, where conservation is observed up to $\mathcal{O}(10^{-12})$. The average rank over all spatial nodes $x_i$ remains very low over time, as also shown in Figure \ref{fig:vfp-conservation-and-rank}.

The second-order temporal accuracy of the IMEX222 method used for the time integration is observed in Figure \ref{fig:vfp-vperpint-temporal-error}. The $L^1$ errors are shown for $f$, $n$, $u_{||}$, and $T_{\alpha}$. Lastly, the expected second-order accuracy in $\vperp$ and $\vpar$ is observed in Tables \ref{tab:vfp-vperp-errors}-\ref{tab:vfp-vpara-errors}, where the velocity mesh is fixed to $N_v=200$ in one variable, and the other variable refines the mesh with $N_v=48,80,144,240$. These mesh values were chosen so that the cell centers aligned with the grid of the reference solution, which used $N_v=720$ in the variable of interest.

\begin{table}[h!]
\centering
\begin{tabular}{|c|cc|cc|}
\hline
$N_{\perp} \times N_{||}$ & $L^1$ Error & Order & $L^{\infty}$ Error & Order \\
\hline
$48  \times 200$ & 1.25E-06 & --   & 1.71E-04 & --   \\
$80  \times 200$ & 4.10E-07 & 2.18 & 5.32E-05 & 2.29 \\
$144 \times 200$ & 1.18E-07 & 2.13 & 1.46E-05 & 2.20 \\
$240 \times 200$ & 3.80E-08 & 2.19 & 4.60E-06 & 2.25 \\
\hline
\end{tabular}
\caption{Error table for the Mach 5 steady-state shock under mesh refinement in $\vperp$. Mesh $N_x=160$, $N_{||}=200$, time-step $\Delta t=0.2$, truncation tolerance $\vartheta=10^{-8}$, and final time $T_f=1$.}
\label{tab:vfp-vperp-errors}
\end{table}

\begin{table}[h!]
\centering
\begin{tabular}{|c|cc|cc|}
\hline
$N_{\perp} \times N_{||}$ & $L^1$ Error & Order & $L^{\infty}$ Error & Order \\
\hline
$200  \times 48$ & 1.69E-05 & --   & 8.90E-03 & --   \\
$200  \times 80$ & 5.98E-06 & 2.04 & 3.20E-03 & 1.98 \\
$200 \times 144$ & 1.79E-06 & 2.06 & 1.00E-03 & 2.04 \\
$200 \times 240$ & 5.90E-07 & 2.15 & 3.00E-04 & 2.15 \\
\hline
\end{tabular}
\caption{Error table for the Mach 5 steady-state shock under mesh refinement in $\vpar$. Mesh $N_x=160$, $N_{\perp}=200$, time-step $\Delta t=0.2$, truncation tolerance $\vartheta=10^{-8}$, and final time $T_f=1$.}
\label{tab:vfp-vpara-errors}
\end{table}

\subsection{Weak Landau damping}

\begin{figure}[h!]
\centering
    \begin{minipage}{0.49\linewidth}
        \centering
        \includegraphics[width=\linewidth]{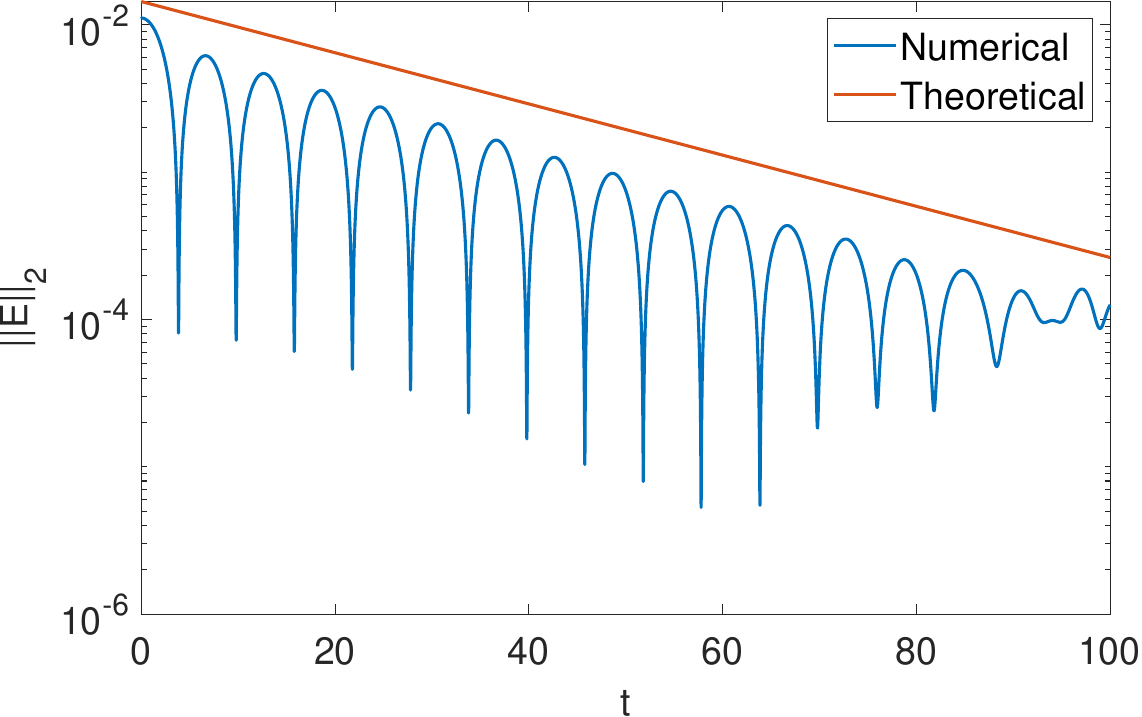}  
    \end{minipage}
    \hfill
    \begin{minipage}{0.49\linewidth}
        \centering
        \includegraphics[width=0.95\linewidth]{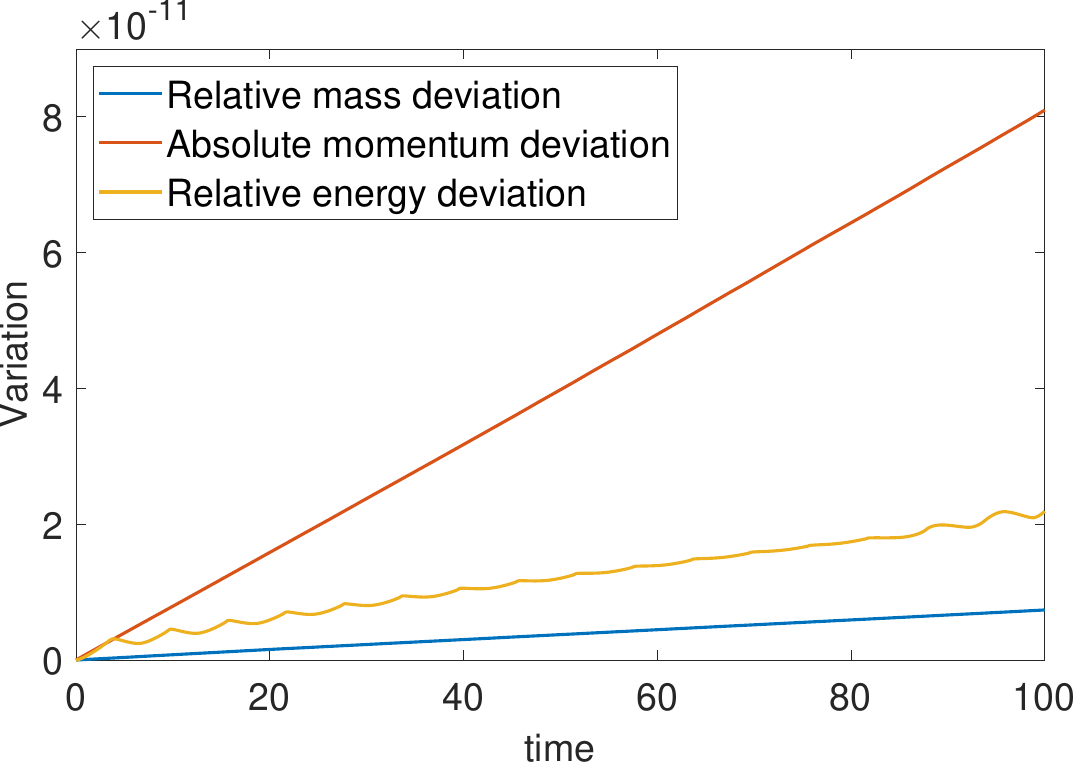}
    \end{minipage}
    \caption{The $L^2$ norm of the electric field \text{(left)}, and the change in mass, momentum, and energy \textit{(right)}. Mesh $N_x=300$, $N_{\perp}=64$, $N_{||}=128$, time-step $\Delta t=0.005$, final time $T_f=100$, and truncation tolerance $\vartheta=10^{-8}$.}
    \label{fig:landau}
\end{figure}

We test the robustness of our method on a collisionless weak Landau damping problem for an ion acoustic wave with kinetic ions and near-isothermal (high-conductivity) electrons, for which a hybrid model is appropriate. We consider equations \eqref{eq:VFP1d2v}-\eqref{eq:fluidelectron} setting $C_{\alpha\alpha}=C_{\alpha e}=0$, $W_{e\alpha}=0$, and scale the thermal conductivity $\kappa_{||}$ by a factor of 10. We approximate the isothermal-electron limit by inflating $\kappa_{||}$ so that $T_e$ remains nearly uniform. The domain is $(x,\vperp,\vpar)\in(0,16)\times(0,6)\times(-6,6)$ with mesh $N_x=300$, $N_{\perp}=64$, $N_{||}=128$. The proton-electron plasma is assumed to initially have spatially homogeneous temperature distributions $T_{\alpha}(x,t=0)=0.2$ and $T_e(x,t=0)=1$, zero drift $u_{||}=0$, and an acoustic wave generated by an ion density perturbation $n(x,t=0)=1+0.01\sin{(2\pi x/16)}$. We plot the $L^2$ norm of the electric field $\norm{E_{||}(x,t)}_2$ and the change in mass, momentum, and energy in Figure \ref{fig:landau}. We ran the simulation to time $T_f=100$ using IMEX222, time-step $\Delta t=0.005$, and truncation tolerance $\vartheta=10^{-8}$. The QCM procedure is unnecessary since the plasma is collisionless. As seen in Figure \ref{fig:landau}, there is good agreement between the simulation and theoretical ion Landau damping rate of $\gamma_{LD}=-0.04$ \cite{taitano_2018_jcp_r_adaptivity,ValentiniTCHM2007hybrid}. In addition, mass, momentum, and energy are conserved up to $\mathcal{O}(10^{-11})$. Since the Newton solver for the fluid equations converges slower for the Vlasov equation than the VFP equation, we made the Newton tolerance $\norm{\mathbf{R}_{\ell}}_{\infty}\leq 10^{-12}$ to maintain efficiency while remaining very conservative; decreasing the tolerance slightly improves the conservation.

\section{Conclusion}

In this paper, we presented a rank-adaptive framework for solving the 1D-2V hybrid Vlasov-Fokker-Planck model. At each spatial node, an implicit-explicit rank-adaptive solver can be used in velocity space to update the kinetic-ion distribution function. We used the recently developed Reduced Augmentation Implicit Low-rank (RAIL) method as our choice of solver, but other suitable solvers can be used just as easily. The Local Macroscopic Conservative (LoMaC) method is used as a post-processing step to truncate the distribution function while conserving mass, momentum, and energy. The extension from Cartesian to azimuthally symmetric cylindrical coordinates was shown for both the RAIL and LoMaC methods. The structure-preserving Chang-Cooper discretization of the LBFP operator was also used. Numerical tests verified the order of accuracy, conservation, and robustness of the method, particularly for long-time simulations of the standing shock and weak Landau damping problems. The flexibility in the choice of implicit-explicit rank-adaptive solver offers an avenue for extending our proposed method to higher-dimensional hybrid kinetic models.

Nonlinear high resolution methods are nontrivial in a majority of rank-adaptive frameworks since one needs to work with factorized solutions over neighboring cells. Minmod and WENO-type limiters have been applied in other rank-adaptive methods for similar problems, such as the discontinuous Galerkin framework \cite{GalindoNPQT2025nodal} and sketching-based low-rank framework \cite{SandsQHZ2025adaptive,ZhengHCQ2025semi,ZhengSHCQ2025semi}. Applying such limiters in projection-based methods such as the RAIL framework and other dynamical low-rank type methods is still a relatively open area of research. Our future work also includes incorporating similar ideas to those in \cite{SandsQHZ2025adaptive,ZhengHCQ2025semi,ZhengSHCQ2025semi} to the current method. Lastly, as with most numerical methods that work over fixed structured grids, the time-stepping size is restricted by the CFL condition. Modifications to a semi-Lagrangian or Eulerian-Lagrangian framework could alleviate this time-step restriction, and is another point of future direction.

\subsection{Data Availability} The source code generated and analyzed in the current study is available at \url{https://github.com/dylan-jacobs/HybridVFPSolver}.

\subsection{Competing Interests} The authors have no competing interests.

\subsection{Acknowledgments} Joseph Nakao was supported by the Swarthmore College Faculty Research Support Grant. William Taitano was supported by the DOE Office of Science, Advanced Scientific Computing Research (ASCR) through the Mathematical Multifaceted Integrated Capability Centers program. William Taitano also acknowledges Triad National Security, LLC contract 89233218CNA000001. Early discussions that led to this work were made while Joseph Nakao was a visiting researcher at the Air Force Research Laboratory in 2022 as a part of the Air Force Summer Faculty Fellowship Program, under the co-mentorship of Alexander Alekseenko (California State University at Northridge) and William Taitano, who at the time was at the Air Force Research Laboratory.



\medskip

\printbibliography

@article{xu2010unified,
  title={A unified gas-kinetic scheme for continuum and rarefied flows},
  author={Xu, Kun and Huang, Juan-Chen},
  journal={Journal of Computational Physics},
  volume={229},
  number={20},
  pages={7747--7764},
  year={2010},
  publisher={Elsevier}
}

@article{Chen2015,
  title={A comparative study of an asymptotic preserving scheme and unified gas-kinetic scheme in continuum flow limit},
  author={Chen, Songze and Xu, Kun},
  journal={Journal of Computational Physics},
  volume={288},
  pages={52--65},
  year={2015},
  month={May}
}

@article{Wei2026,
  title={Unified gas-kinetic scheme for reactive flow with multi-scale transport and chemical non-equilibrium},
  author={Wei, Yufeng and Cao, Junzhe and Xu, Kun},
  journal={Journal of Computational Physics},
  volume={546},
  pages={114514},
  year={2026},
  month={February}
}

@article{rosenbluth_1957_prl_fp,
  title={Fokker-{P}lanck equation for an inverse-square force},
  author={Rosenbluth, M.N and MacDonald, W.M and Judd, D.L},
  journal={Physical Review},
  volume={107},
  number={1},
  pages={1--6},
  year={1957}
}

@article{johnson_pre_2024_omega,
  title={Impact of mid-{Z} gas fill on dynamics and performance of shock-driven implosions at the {OMEGA} laser},
  author={Gatu Johnson, M and others},
  journal={Physical Review E},
  volume={109},
  number={6},
  pages={065201},
  year={2024}
}

@article{mannion_pre_2023_omega,
  title={Evidence of non-Maxwellian ion velocity distributions in spherical shock-driven implosions},
  author={Mannion, O.M and others},
  journal={Physical Review E},
  volume={108},
  number={3},
  pages={035201},
  year={2023}
}

@article{taitano_2016_jcp_r_adaptivity,
  title={An adaptive, conservative 0{D}-2{V} multispecies {R}osenbluth-{F}okker-{P}lanck solver for arbitrarily disparate mass and temperature regimes},
  author={Taitano, W.T and Chac{\'o}n, L and Simakov, A.N},
  journal={Journal of Computational Physics},
  volume={318},
  number={},
  pages={391--420},
  year={2016}
}

@article{taitano_2018_jcp_r_adaptivity,
  title={An adaptive, implicit, conservative, 1D-2V multi-species Vlasov--Fokker--Planck multi-scale solver in planar geometry},
  author={Taitano, William T and Chac{\'o}n, Luis and Simakov, Andrei N},
  journal={Journal of Computational Physics},
  volume={365},
  pages={173--205},
  year={2018},
  publisher={Elsevier}
}

@article{taitano_2021_jcp_r_adaptivity,
  title={An {E}ulerian {V}lasov-{F}okker-{P}lanck algorithm for spherical implosion simulations of inertial confinement fusion capsules},
  author={Taitano, W and Keenan, B and Chac{\'o}n, L and Anderson, S and Hammer, H and Simakov, A},
  journal={Computer Physics Communications},
  volume={263},
  number={},
  pages={107861},
  year={2021}
}

@article{kolobov_aip_cp_2019_amr_vlasov,
  title={Boltzmann-{F}okker-{P}lanck kinetic solver with adaptive mesh in phase space},
  author={Kolobov, V and Arslanbekov, R and Levko, D},
  journal={AIP Conf. Proc.},
  volume={2132},
  number={},
  pages={060011},
  year={2019}
}

@article{bennoune_jcp_2008_mmd,
  title={Uniformly stable numerical schemes for the Boltzmann equation preserving the compressible Navier-Stokes asymptotics},
  author={Bennoune, M and Lemou, M and Mieussens, L},
  journal={Journal of Computational Physics},
  volume={227},
  number={},
  pages={3781--3803},
  year={2008}
}

@article{gamba_jcp_2019_mmd,
  title={Micro-macro decomposition based asymptotic-preserving numerical schemes and numerical moments conservation for collisional nonlinear kinetic equations},
  author={Gamba, I and Jin, S and Liu, L},
  journal={Journal of Computational Physics},
  volume={382},
  number={},
  pages={264--290},
  year={2019}
}

@article{lemou_jsc_2008_mmd,
  title={A new asymptotic preserving scheme based on micro-macro formulation for linear kinetic equations in the diffusion limit},
  author={Lemou, M and Mieussens, L},
  journal={SIAM J. Sci. Comput.},
  volume={31},
  number={1},
  pages={334--368},
  year={2008}
}

@article{liu_cmp_2004_mmd,
  title={Micro-macro decompositions and positivity of shock profiles},
  author={Liu, T and Yu, S},
  journal={Comm. Math. Phys.},
  volume={246},
  number={},
  pages={133--179},
  year={2004}
}

@article{taitano_jpp_2026_ck,
  title={A novel conditional formulation of the Vlasov-Ampere equations: a conservative, positivity, asymptotic, and Gauss law preserving scheme},
  author={Taitano, W and Burby, J and Alekseenko, A},
  journal={Journal of Plasma Physics},
  volume={92},
  number={},
  pages={E45},
  year={2026}
}

@article{holo_review_2017,
  title={Multiscale high-order/low-order (HOLO) algorithms and applications},
  author={Chac{\'o}n, L and Chen, G and Knoll, D and Newman, C and Park, H and Taitano, W and Willert, J and Womeldorff, G},
  journal={Journal of Computational Physics},
  volume={330},
  number={},
  pages={21--45},
  year={2017}
}

@article{taitano_vfp_holo_jcp_2015,
  title={Charge-and-energy conserving moment-based accelerator for multi-species Vlasov-Fokker-Planck-Amp{`e}re system, part II: collisional aspects},
  author={Taitano, W and Chac{\'o}n, L and Knoll, D},
  journal={Journal of Computational Physics},
  volume={284},
  number={},
  pages={737--757},
  year={2015}
}

@article{taitano_bgk_holo_jctt_2014,
  title={Moment-based acceleration for neutral gas kinetics with {BGK} collision operator},
  author={Taitano, W.T and Knoll, D.A and Chac{\'o}n, L and Reisner, J and Prinja, A},
  journal={J. Comput. Theor. Transp.},
  volume={43},
  number={1-7},
  pages={83--108},
  year={2014}
}

@article{willert_neutraon_holo_jcp_2014,
  title={Leveraging Anderson acceleration for improved convergence of iterative solutions to transport systems},
  author={Willert, J and Taitano, W.T and Knoll, D},
  journal={Journal of Computational Physics},
  volume={273},
  number={},
  pages={278--286},
  year={2014}
}

@article{LenardB1958,
  title={Plasma oscillations with diffusion in velocity space},
  author={Lenard, Andrew and Bernstein, Ira B},
  journal={Physical Review},
  volume={112},
  number={5},
  pages={1456},
  year={1958},
  publisher={APS}
}

@article{Dougherty1964,
  title={Model Fokker-Planck equation for a plasma and its solution},
  author={Dougherty, JP},
  journal={The Physics of Fluids},
  volume={7},
  number={11},
  pages={1788--1799},
  year={1964},
  publisher={AIP Publishing}
}

@article{ChangC1970practical,
  title={A practical difference scheme for Fokker-Planck equations},
  author={Chang, JS and Cooper, G},
  journal={Journal of Computational Physics},
  volume={6},
  number={1},
  pages={1--16},
  year={1970},
  publisher={Elsevier}
}

@article{PareschiZ2018structure,
  title={Structure preserving schemes for nonlinear Fokker--Planck equations and applications},
  author={Pareschi, Lorenzo and Zanella, Mattia},
  journal={Journal of Scientific Computing},
  volume={74},
  number={3},
  pages={1575--1600},
  year={2018},
  publisher={Springer}
}

@article{KahzaCTQH2026structure,
  title={A Structure-Preserving Penalization Method for the Single-species Rosenbluth-Fokker-Planck Equation},
  author={Kahza, Hamad El and Chac{\'o}n, Luis and Taitano, William and Qiu, Jingmei and Hu, Jingwei},
  journal={arXiv preprint arXiv:2601.08006},
  year={2026}
}

@article{VidalMCL1993ion,
  title={Ion kinetic simulations of the formation and propagation of a planar collisional shock wave in a plasma},
  author={Vidal, F and Matte, JP and Casanova, M and Larroche, O},
  journal={Physics of Fluids B: Plasma Physics},
  volume={5},
  number={9},
  pages={3182--3190},
  year={1993},
  publisher={American Institute of Physics}
}

@article{TaitanoCS2021conservative,
  title={A conservative phase-space moving-grid strategy for a 1D-2V Vlasov--Fokker--Planck Solver},
  author={Taitano, William T and Chacon, Luis and Simakov, Andrei N and Anderson, Steven E},
  journal={Computer Physics Communications},
  volume={258},
  pages={107547},
  year={2021},
  publisher={Elsevier}
}

@article{taitano2015mass,
  title={A mass, momentum, and energy conserving, fully implicit, scalable algorithm for the multi-dimensional, multi-species Rosenbluth--Fokker--Planck equation},
  author={Taitano, William T and Chac{\'o}n, Luis and Simakov, AN and Molvig, K},
  journal={Journal of Computational Physics},
  volume={297},
  pages={357--380},
  year={2015},
  publisher={Elsevier}
}

@article{TaitanoCS2017equilibrium,
  title={An equilibrium-preserving discretization for the nonlinear Rosenbluth--Fokker--Planck operator in arbitrary multi-dimensional geometry},
  author={Taitano, William T and Chacon, Luis and Simakov, Andrei Nikolaevich},
  journal={Journal of Computational Physics},
  volume={339},
  pages={453--460},
  year={2017},
  publisher={Elsevier}
}

@article{GuoQ2022low,
  title={A low rank tensor representation of linear transport and nonlinear Vlasov solutions and their associated flow maps},
  author={Guo, Wei and Qiu, Jing-Mei},
  journal={Journal of Computational Physics},
  volume={458},
  pages={111089},
  year={2022},
  publisher={Elsevier}
}

@article{SimakovM2014electron,
  title={Electron transport in a collisional plasma with multiple ion species},
  author={Simakov, Andrei N and Molvig, Kim},
  journal={Physics of Plasmas},
  volume={21},
  number={2},
  year={2014},
  publisher={AIP Publishing}
}

@book{HazeltineM2003plasma,
  title={Plasma confinement},
  author={Hazeltine, Richard D and Meiss, James D},
  year={2003},
  publisher={Courier Corporation}
}

@article{ValentiniTCHM2007hybrid,
  title={A hybrid-Vlasov model based on the current advance method for the simulation of collisionless magnetized plasma},
  author={Valentini, Francesco and Tr{\'a}vn{\'\i}{\v{c}}ek, P and Califano, Francesco and Hellinger, Petr and Mangeney, Andr{\'e}},
  journal={Journal of Computational Physics},
  volume={225},
  number={1},
  pages={753--770},
  year={2007},
  publisher={Elsevier}
}

@article{GalindoNPQT2025nodal,
  title={A Nodal Discontinuous Galerkin Method with Low-Rank Velocity Space Representation for the Multi-Scale BGK Model},
  author={Galindo-Olarte, Andres and Nakao, Joseph and Pasha, Mirjeta and Qiu, Jing-Mei and Taitano, William},
  journal={arXiv preprint arXiv:2508.16564},
  year={2025}
}

@article{NakaoQE2025reduced,
  title={Reduced Augmentation Implicit Low-rank (RAIL) integrators for advection-diffusion and Fokker--Planck models},
  author={Nakao, Joseph and Qiu, Jing-Mei and Einkemmer, Lukas},
  journal={SIAM Journal on Scientific Computing},
  volume={47},
  number={2},
  pages={A1145--A1169},
  year={2025},
  publisher={SIAM}
}

@article{NakaoCE2026low,
  title={A low-rank, high-order implicit-explicit integrator for three-dimensional convection-diffusion equations},
  author={Nakao, Joseph and Ceruti, Gianluca and Einkemmer, Lukas},
  journal={Computer Physics Communications},
  volume={325},
  pages={110163},
  year={2026},
  publisher={Elsevier}
}

@article{MengAC2025preconditioning,
  title={Preconditioning low rank generalized minimal residual method (gmres) for implicit discretizations of matrix differential equations},
  author={Meng, Shixu and Appel{\"o}, Daniel and Cheng, Yingda},
  journal={SIAM Journal on Matrix Analysis and Applications},
  volume={46},
  number={4},
  pages={2475--2500},
  year={2025},
  publisher={SIAM}
}

@article{KahzaTQC2024krylov,
  title={Krylov-based adaptive-rank implicit time integrators for stiff problems with application to nonlinear Fokker-Planck kinetic models},
  author={El Kahza, Hamad and Taitano, William and Qiu, Jing-Mei and Chac{\'o}n, Luis},
  journal={Journal of Computational Physics},
  volume={518},
  pages={113332},
  year={2024},
  publisher={Elsevier}
}

@article{CerutiEKL2024robust,
  title={A robust second-order low-rank BUG integrator based on the midpoint rule},
  author={Ceruti, Gianluca and Einkemmer, Lukas and Kusch, Jonas and Lubich, Christian},
  journal={BIT Numerical Mathematics},
  volume={64},
  number={3},
  pages={30},
  year={2024},
  publisher={Springer}
}

@article{LiJ2026high,
  title={High-Order Implicit Low-Rank Method with Spectral Deferred Correction for Matrix Differential Equations},
  author={Li, Shun and Jiang, Yan},
  journal={Journal of Scientific Computing},
  volume={106},
  number={3},
  pages={75},
  year={2026},
  publisher={Springer}
}

@article{GuoQ2024local,
  title={A Local Macroscopic Conservative (LoMaC) low rank tensor method for the Vlasov dynamics},
  author={Guo, Wei and Qiu, Jing-Mei},
  journal={Journal of Scientific Computing},
  volume={101},
  number={3},
  pages={61},
  year={2024},
  publisher={Springer}
}

@article{SandsQHZ2025adaptive,
  title={An adaptive-rank approach with greedy sampling for multi-scale BGK equations},
  author={Sands, William A and Qiu, Jing-Mei and Hayes, Daniel and Zheng, Nanyi},
  journal={Journal of Computational Physics},
  pages={114523},
  year={2025},
  publisher={Elsevier}
}

@article{ZhengHCQ2025semi,
  title={A semi-Lagrangian adaptive-rank (SLAR) method for linear advection and nonlinear Vlasov-Poisson system},
  author={Zheng, Nanyi and Hayes, Daniel and Christlieb, Andrew and Qiu, Jing-Mei},
  journal={Journal of Computational Physics},
  volume={532},
  pages={113970},
  year={2025},
  publisher={Elsevier}
}

@article{ZhengSHCQ2025semi,
  title={A Semi-Lagrangian Adaptive Rank (SLAR) Method for High-Dimensional Vlasov Dynamics},
  author={Zheng, Nanyi and Sands, William A and Hayes, Daniel and Christlieb, Andrew J and Qiu, Jing-Mei},
  journal={arXiv preprint arXiv:2510.24861},
  year={2025}
}

@article{AscherRS1997implicit,
  title={Implicit-explicit Runge-Kutta methods for time-dependent partial differential equations},
  author={Ascher, Uri M and Ruuth, Steven J and Spiteri, Raymond J},
  journal={Applied Numerical Mathematics},
  volume={25},
  number={2-3},
  pages={151--167},
  year={1997},
  publisher={Elsevier}
}

@article{Simoncini2016computational,
  title={Computational methods for linear matrix equations},
  author={Simoncini, Valeria},
  journal={SIAM Review},
  volume={58},
  number={3},
  pages={377--441},
  year={2016},
  publisher={SIAM}
}

@article{hackbusch2009new,
	title={A new scheme for the tensor representation},
	author={Hackbusch, Wolfgang and K{\"u}hn, Stefan},
	journal={Journal of Fourier analysis and applications},
	volume={15},
	number={5},
	pages={706--722},
	year={2009},
	publisher={Springer}
}

@article{grasedyck2013literature,
	title={A literature survey of low-rank tensor approximation techniques},
	author={Grasedyck, Lars and Kressner, Daniel and Tobler, Christine},
	journal={GAMM-Mitteilungen},
	volume={36},
	number={1},
	pages={53--78},
	year={2013},
	publisher={Wiley Online Library}
}

@article{kolda2009tensor,
  title={Tensor decompositions and applications},
  author={Kolda, Tamara G and Bader, Brett W},
  journal={SIAM review},
  volume={51},
  number={3},
  pages={455--500},
  year={2009},
  publisher={SIAM}
}

@article{kormann2017low,
  title={Low-rank tensor discretization for high-dimensional problems},
  author={Kormann, Katharina},
  year={2017},
  journal={Vorlesung (SS 2017)}
}

@article{bachmayr2023low,
  title={Low-rank tensor methods for partial differential equations},
  author={Bachmayr, Markus},
  journal={Acta Numerica},
  volume={32},
  pages={1--121},
  year={2023},
  publisher={Cambridge University Press}
}

@article{einkemmer2025review,
  title={A review of low-rank methods for time-dependent kinetic simulations},
  author={Einkemmer, Lukas and Kormann, Katharina and Kusch, Jonas and McClarren, Ryan G and Qiu, Jing-Mei},
  journal={Journal of Computational Physics},
  pages={114191},
  year={2025},
  publisher={Elsevier}
}

@article{nelson2026local,
  title={A Local Macroscopic Conservative Low-Rank Discontinuous Galerkin Method for the Vlasov-Poisson Equation with Dougherty-Fokker-Planck Collisions},
  author={Nelson, Austin and Guo, Wei and Guthrey, Pierson},
  journal={arXiv preprint arXiv:2607.19264},
  year={2026}
}

@article{carrel2026sketch,
  title={Sketch low-rank dynamics: orthogonal vs. oblique projections},
  author={Carrel, Benjamin and Grigori, Laura},
  journal={arXiv preprint arXiv:2607.03402},
  year={2026}
}

@article{dektor2021rank,
  title={Rank-adaptive tensor methods for high-dimensional nonlinear PDEs},
  author={Dektor, Alec and Rodgers, Abram and Venturi, Daniele},
  journal={Journal of Scientific Computing},
  volume={88},
  number={2},
  pages={36},
  year={2021},
  publisher={Springer}
}

@article{wang2026implicit,
  title={Implicit Dynamical Tensor Train Approximation for Kinetic Equations with Stiff Fokker--Planck Collisions},
  author={Wang, Geshuo and Hu, Jingwei},
  journal={arXiv preprint arXiv:2605.15382},
  year={2026}
}

@article{wang2026dynamical,
  title={Dynamical tensor train approximation for kinetic equations},
  author={Wang, Geshuo and Hu, Jingwei},
  journal={Journal of Computational Physics},
  pages={114884},
  year={2026},
  publisher={Elsevier}
}

@article{zhang2026separable,
  title={A Separable and Asymptotic-Preserving Dynamical Low-Rank Method for the Vlasov--Poisson--Fokker--Planck System},
  author={Zhang, Shiheng and Hu, Jingwei},
  journal={arXiv preprint arXiv:2601.11900},
  year={2026}
}

@article{coughlin2024robust,
  title={Robust and conservative dynamical low-rank methods for the Vlasov equation via a novel macro-micro decomposition},
  author={Coughlin, Jack and Hu, Jingwei and Shumlak, Uri},
  journal={Journal of Computational Physics},
  volume={509},
  pages={113055},
  year={2024},
  publisher={Elsevier}
}

@article{coughlin2022efficient,
  title={Efficient dynamical low-rank approximation for the Vlasov-Amp{\`e}re-Fokker-Planck system},
  author={Coughlin, Jack and Hu, Jingwei},
  journal={Journal of Computational Physics},
  volume={470},
  pages={111590},
  year={2022},
  publisher={Elsevier}
}

@article{chertkov2021solution,
  title={Solution of the Fokker--Planck equation by cross approximation method in the tensor train format},
  author={Chertkov, Andrei and Oseledets, Ivan},
  journal={Frontiers in Artificial Intelligence},
  volume={4},
  pages={668215},
  year={2021},
  publisher={Frontiers Media SA}
}

@article{dolgov2012fast,
  title={Fast solution of parabolic problems in the tensor train/quantized tensor train format with initial application to the Fokker--Planck equation},
  author={Dolgov, Sergey V and Khoromskij, Boris N and Oseledets, Ivan V},
  journal={SIAM Journal on Scientific Computing},
  volume={34},
  number={6},
  pages={A3016--A3038},
  year={2012},
  publisher={SIAM}
}

\end{document}